\documentclass{amsart}
\usepackage[toc, page]{appendix}
\usepackage{comment}
\usepackage{amsmath, amsthm, amssymb, color}
\usepackage{graphicx}
\usepackage{mwe}
\usepackage{stmaryrd}
\usepackage{a4wide,bm}
\usepackage{amssymb}
\usepackage{euscript}
\usepackage{xcolor}
\usepackage{url}
\usepackage{stmaryrd}
\usepackage{microtype}
\usepackage{graphicx}
\usepackage{amsmath}
\usepackage{subcaption}
\usepackage{float}
\usepackage{enumitem}   
\usepackage{mathrsfs}
\usepackage{diagbox}
\usepackage[all]{xy}

\usepackage[T1]{fontenc}
\usepackage{textcomp}
\usepackage{tikz}
\usepackage{capt-of}
\usepackage{lineno}
\usepackage{indentfirst}
\usepackage{mathtools}
\usepackage{hyperref}

\hypersetup{colorlinks=true,
	linkcolor=red,
	filecolor=blue,  
	citecolor=blue, 
	urlcolor=blue}

\usepackage[capitalize,noabbrev]{cleveref}

\newcommand{\circst}{\varoast}

\newcommand{\R}{\mathbb{R}}

\newcommand{\Z}{\mathbb{Z}}

\newcommand{\E}{\mathbb{E}}
\newcommand{\T}{\mathbb{T}}

\renewcommand{\P}{\mathbb{P}}

\newcommand{\Lcal}{\mathcal{L}}

\newcommand{\norm}[1]{\left\lVert#1\right\rVert}

\newcommand*\varhrulefill[1][2pt]{\leavevmode\leaders\hrule height#1\hfill\kern0pt}
\newcommand*{\Scale}[2][4]{\scalebox{#1}{$#2$}}%

\newtheorem{theorem}{Theorem}[section]
\newtheorem{lemma}[theorem]{Lemma}
\newtheorem{corollary}[theorem]{Corollary}

\theoremstyle{definition}

\newtheorem{example}[theorem]{Example}

\newtheorem{proposition}[theorem]{Proposition}

\newtheorem{conjecture}[theorem]{Conjecture}
\newtheorem{question}[theorem]{Question}
\theoremstyle{remark}
\newtheorem{remark}[theorem]{Remark}

\numberwithin{equation}{section}

\begin{document}
	
	\title[The Bernoulli clock]{The Bernoulli clock: probabilistic and combinatorial interpretations of the Bernoulli polynomials by circular convolution}
	
	\author{Yassine El Maazouz}
	\address{Yassine El Maazouz, U.C. Berkeley, Department of statistics, 335 Evans Hall \#3860 Berkeley, CA 94720 U.S.A. }
	\email{yassine.el-maazouz@berkeley.edu}

	\author{Jim Pitman}
	\address{Jim Pitman, U.C. Berkeley, Department of statistics, 367 Evans Hall \#3860 Berkeley, CA 94720 U.S.A.}
	\email{pitman@berkeley.edu}
	
	
	\date{\today}

	\begin{abstract}
		The factorially normalized Bernoulli polynomials $b_n(x) = B_n(x)/n!$ are known to be characterized by  $b_0(x) = 1$ and $b_n(x)$ for $n >0$ is the anti-derivative of $b_{n-1}(x)$ subject to $\int_0^1 b_n(x) dx = 0$.
		We offer a related characterization: $b_1(x) = x - 1/2$ and $(-1)^{n-1} b_n(x)$ for $n >0$ is the $n$-fold circular convolution of $b_1(x)$ with itself.
		Equivalently, $1 - 2^n b_n(x)$ is the probability density at $x \in (0,1)$ of the fractional part of a sum of $n$ independent random variables, each with the beta$(1,2)$ probability density $2(1-x)$ at $x \in (0,1)$.
		This result has a novel combinatorial analog, the {\em Bernoulli clock}: mark the hours of a $2 n$ hour clock by a uniform random permutation of the multiset $\{1,1, 2,2, \ldots, n,n\}$, meaning pick two different hours uniformly at random from the $2 n$ hours and mark them $1$, then pick two different hours uniformly at random from the remaining $2 n - 2$ hours and mark them $2$, and so on.
		Starting from hour $0 = 2n$, move clockwise to the first hour marked $1$, continue clockwise to the first hour marked $2$, and so on, continuing clockwise around the Bernoulli clock until the first of the two hours marked $n$ is encountered, at a random hour $I_n$ between $1$ and $2n$. 
		We show that for each positive integer $n$,  the event $( I_n = 1)$ has probability $(1 - 2^n b_n(0))/(2n)$, where $n! b_n(0) = B_n(0)$ is the $n$th Bernoulli number.
		For $ 1 \le k \le 2 n$, the difference $\delta_n(k):= 1/(2n) - \P( I_n = k)$ is a polynomial function of $k$ with the surprising symmetry $\delta_n( 2 n + 1 - k)  = (-1)^n \delta_n(k)$, which is a combinatorial analog of the well known symmetry of Bernoulli polynomials $b_n(1-x) = (-1)^n b_n(x)$.
	\end{abstract}
	
	\maketitle
	
	\setcounter{tocdepth}{1}
	\tableofcontents
	
	\section{Introduction}
	
	The \emph{Bernoulli polynomials} $(B_n(x))_{n \geq 0}$ are a special sequence of univariate polynomials with rational coefficients. They are named after the Swiss mathematician Jakob Bernoulli (1654–1705), who (in his \textit{Ars Conjectandi} published posthumously in Basel 1713) found the sum of $m$th powers of the first $n$ positive integers using the instance $x = 1$ of the {\em power sum formula}
	\begin{equation}\label{eq:powerSum}
		\sum_{k=0}^{n-1} (x+k)^{m} = \frac{ B_{m+1} (x+n) - B_{m+1}(x) }{ m + 1} , \qquad (n = 1,2, \ldots, m = 0,1,2, \ldots) .
	\end{equation}
	The evaluations $B_m:= B_m(0)$ and $B_m(1) = (-1)^{m} B_m$  are known as the \emph{Bernoulli numbers}, from which the polynomials are recovered as
	\begin{equation}\label{eq:bernexp}
		B_n(x) = \sum_{k=0}^{n} \binom{n}{k} B_{n-k} \ x^k.
	\end{equation}
	These polynomials have been well studied, starting from the early work of Faulhaber, Bernoulli, Seki and Euler in the $17$th and early $18$th centuries. They can be defined in multiple ways. For example, Euler defined the Bernoulli polynomials by their \emph{exponential generating function}
	
	\begin{equation} \label{eq:Euler_Gen_Func}
		B(x,\lambda) \coloneqq \frac{\lambda e^{\lambda x } }{e^{\lambda} - 1}  = \sum_{n = 0}^{\infty} \frac{ B_n(x) }{n!} \lambda^n \qquad \qquad (|\lambda| < 2 \pi ).
	\end{equation} 
	
	Beyond evaluating power sums, the Bernoulli numbers and polynomials are useful in other contexts and appear in many areas in mathematics, among which we mention number theory \cite{BernoulliZetaFunctions,EulerAndZeta74,Mazur11,agoh_fermats_1982}, Lie theory \cite{BourbakiLieGroupsChapII,RepsOfSU3,GradedLieAlgebras,MagnusExpansion}, algebraic geometry and topology \cite{HirzebruchAlgTopAlgGeom, MilnorKervaire}, probability \cite{LevyPlanarBM,LevyBernoulliNumbers, BernoulliEulerPolysStochasticArea1, BernoulliEulerPolysStochasticArea2,sun_moment_2007,pitman_infinitely_2003,biane_probability_2001} and numerical approximation \cite{Steffensen,PhillipsApprox}.
	
	The {\em factorially normalized Bernoulli polynomials} $b_n(x):= B_n(x)/n!$
	can also be defined inductively as follows (see \cite[\S 9.5]{MontgomeryEarlyFourierAnalysis}). Beginning with $b_0(x) = B_0(x) = 1$, 
	for each positive integer $n$,
	the function $x \mapsto b_n(x)$ is the unique antiderivative of $x \mapsto b_{n-1}(x)$ that integrates to $0$ over $[0,1]$:
	
	\begin{equation} \label{eq:intchar}
		b_0(x) = 1,  \quad \frac{ d } { d x } b_n(x) = b_{n-1} (x) \quad \mbox{and} \quad \int_0^1 b_n(x) = 0 \qquad (n > 0) .
	\end{equation}
	So the first few polynomials $b_n(x)$ are 
	\begin{align*}
		b_0(x) &= 1,  										  &b_1(x) &= x - 1/2,\\
		b_2(x) &= \frac{1}{2!}(x^2 - x - 1/6), 				&b_3(x) &= \frac{1}{3!}(x^3 - 3 x^2/2 + x/2).
	\end{align*}
	As shown in \cite[Theorem 9.7]{MontgomeryEarlyFourierAnalysis} starting from \eqref{eq:intchar}, the functions $f(x) = b_n(x)$ 
	with argument $x \in [0,1)$ are also characterized by the simple form of their {\em Fourier transform}
	\begin{equation}
		\label{eq:ft}
		\widehat{f}(k) \coloneqq \int_{0}^{1} f(x) e^{- 2 \pi i k x} dx  \qquad  (k \in \Z)
	\end{equation}
	which is given by
	\begin{equation} \label{eq:bft}
		\begin{alignedat}{2}
			\widehat{b}_0(k)  &= 1[ k = 0], \qquad   &\mbox{for } k \in \Z ;      \\
			\widehat{b}_n(0) &= 0 \quad \text{and} \quad  \widehat{b}_n(k) = - \frac{1}{(2 \pi i k)^n},   \qquad  &\mbox{for } n > 0  \mbox{ and } k \ne 0,
		\end{alignedat}
	\end{equation}
	with the notation $1[\cdots]$ equal to $1$ if $[\cdots]$ holds and $0$ otherwise. It follows from the Fourier expansion of $b_{n}(x)$:
	\[
	b_n(x) = - \frac{2}{(2\pi)^n} \sum_{k = 1}^{\infty} \frac{1}{k^n} \cos\left(2k\pi x - \frac{n\pi}{2}\right) 
	\]
	that there exists a constant $C > 0$ such that
	\begin{equation} \label{eq:SineFunctionApprox}
		\sup\limits_{0 \leq x \leq 1} \left| (2\pi)^n b_n(x) + 2\cos\left(2 \pi x - \frac{n\pi}{2}\right) \right | \leq C 2^{-n} \quad \text{for } n \geq 2,
	\end{equation}
	see \cite{lehmer_maxima_1940}. So as $n \uparrow \infty$ the polynomials $b_n(x)$ looks like shifted cosine functions. Besides \eqref{eq:Euler_Gen_Func} and \eqref{eq:intchar}, several other characterizations of the Bernoulli polynomials are described in \cite{lehmer_new_1988,CDG06}.
	
	\medskip
	
	{
		This article first draws attention to a simple characterization of the Bernoulli polynomials by {\em circular convolution} and, more importantly, provides an interesting probabilistic and combinatorial interpretation in terms of statistics of random permutations of a multiset.
	}
	
	For a pair of functions $f = f(u)$ and $g = g(u)$, defined for $u$ in $[0,1)$ identified with the circle group $\T:= \R/\Z = [0,1)$, 
	with $f$ and $g$ integrable with respect to Lebesgue measure on $\T$, their {\em circular convolution} $f\circst g$ is the function
	\begin{equation} \label{eq:circonv}
		(f \circst g )(u) = \int_{\T} f(v) g (u-v) dv \qquad  \text{for } u \in \T.
	\end{equation}
	Here $u-v$ is evaluated in the circle group $\T$, that is modulo $1$, and $dv$ is the shift-invariant Lebesgue measure on $\T$ with total measure $1$.
	Iteration of this operation defines the $n$th convolution power $u \mapsto f^{\circst n}(u)$ for each positive integer $n$, each integrable $f$, and $u \in \T$.
	
	\begin{theorem} \label{thm:Bernoulliconv}
		The factorially normalized Bernoulli polynomials $b_n(x) = \frac{B_n(x)}{n!}$ are characterized by:
		\begin{enumerate}[label=(\roman*)]
			\item  $b_0(x) = 1$ and $b_1(x) = x - 1/2$,
			\item  for $n >0$ the $n$-fold circular convolution of $b_1(x)$ with itself is $(-1)^{n-1} b_n(x)$; that is
			\begin{equation} \label{eq:bndef}
				b_n(x) = (-1)^{n-1} b_1^{\circst n}(x).
			\end{equation}
		\end{enumerate} 
	\end{theorem}
	
	In view of the identity $\widehat{f \circst g} = \widehat{f} \  \widehat{g}$,
	\cref{thm:Bernoulliconv} follows from the classical Fourier evaluation \eqref{eq:bft} and uniqueness of the Fourier transform.
	A more elementary proof of \cref{thm:Bernoulliconv}, without Fourier transforms, is provided in Section 
	\ref{Sec:CircularConvOfPoly}. 
	So the Fourier evaluation \eqref{eq:bft} may be regarded as a corollary of \cref{thm:Bernoulliconv}.
	That theorem can also be reformulated as follows:
	\begin{corollary} \label{cor:Bernoulliconv}
		The following identities hold for circular convolution of factorially normalized Bernoulli polynomials:
		\begin{align*}
			b_0(x) \circst b_0(x) &= b_0(x)\\
			b_0(x) \circst b_n(x) &= 0, \quad (n \geq 1),\\
			b_n(x) \circst b_m(x) &= - b_{n + m}(x) ,  \quad (n, m \geq 1).\\
		\end{align*}
	\end{corollary}
	In particular, for positive integers $n$ and $m$, this evaluation of $(b_n \circst b_m)(1)$ yields
	an identity which appears in \cite[p. 31]{norlund_vorlesungen_1924}:
	\begin{equation}
		(-1)^{m} \int_{0}^{1} b_n(u) b_m(u)  du = \int_{0}^{1} b_n(u) b_m(1-u)  du = - b_{n+m}(1) .
	\end{equation}
	Here the first equality is due to the well known {\em reflection symmetry} of the Bernoulli polynomials
	\begin{equation} \label{eq:refl}
		(-1)^m b_m(u) = b_m(1-u) 
		\qquad (m \ge 0 ) 
	\end{equation}
	which is the identity of coefficients of $\lambda^m$ in the elementary identity of Eulerian generating functions
	\begin{equation} \label{eq:eulersym}
		B(u, - \lambda)  = \frac{ (-\lambda) e^{-\lambda u } } { e^{- \lambda} - 1} = \frac{ \lambda e^{\lambda(1-u) } } { e^{\lambda} - 1 } = B(1-u,\lambda).
	\end{equation}

	The rest of this article is organized as follows. \cref{Sec:CircularConvOfPoly} gives an elementary proof for \cref{thm:Bernoulliconv}, and discusses circular convolution of polynomials.
	In \cref{Sec:Prob_Interpretations} we highlight the fact that $1 - 2^n b_n(x)$ is the probability density at $x \in (0,1)$ of the fractional part of a sum of $n$ independent random variables, each with the beta$(1,2)$ probability density $2(1-x)$ at $x \in (0,1)$. Because the minimum of two independent uniform $[0,1]$ variables has this beta$(1,2)$ probability density the circular convolution of $n$ independent beta$(1,2)$ variables is closely related to a continuous model we call the \emph{Bernoulli clock}: Spray the circle $\T = [0,1)$ of circumference $1$ with $2n$ i.i.d uniform positions $U_1, U_1', \dots, U_n, U_n'$ with order statistics
	\[
	U_{1:2n} < \dots < U_{2n:2n}.
	\]
	Starting from the origin $0$, move clockwise to the first of position of the pair $(U_1, U_1')$, continue clockwise to the first position of the pair $(U_2, U_2')$, and so on, continuing clockwise around the circle until the first of the two positions $(U_n,U_n')$ is encountered at a random index $1 \leq I_n \leq 2n$ (i.e. we stop at $U_{I_n:2n}$) after having made a random number $0 \leq D_n \leq n - 1$ turns around the circle. Then for each positive integer $n$,  the event $( I_n = 1)$ has probability 
	$$\P(I_n = 1 ) =  \frac{ 1 - 2^n b_n(0) }{2n}$$ 
	where $n! b_n(0) = B_n(0)$ is the $n$th Bernoulli number. For $ 1 \le k \le 2 n$, the difference 
	$$\delta_{k:2n} := \frac{1}{2n} - \P(I_n = k) $$ 
	is a polynomial function of $k$, which is closely related to $b_n(x)$.
	In particular, this difference has the surprising symmetry $$\delta_{2 n + 1 - k : 2n}  = (-1)^n \delta_{k:2n}, \quad \mbox{for } 1 \leq k \leq 2n$$ which is a combinatorial analog of the reflection symmetry \eqref{eq:refl} for the Bernoulli polynomials. 
	
	Stripping down the clock model, the random variables $I_n$ and $D_n$ are two statistics of permutations of the multiset 
	\begin{equation} \label{eq:ourmulti}
		1^2 \cdots n^2:= \{1,1, 2,2, \ldots, n,n\}.
	\end{equation}
	\cref{Sec:BernoulliClock} discusses the combinatorics behind the distributions of $I_n$ and $D_n$.
	In \cref{Sec:Generalized_Bernoulli_clock} we generalize the Bernoulli clock model to offer a new perspective on the work of Horton and Kurn \cite{HortonKurn} and the more recent work of Clifton et al \cite{CliftonEtAl}. In particular, we provide a probabilistic interpretation for the permutation counting problem in \cite{HortonKurn}  and prove Conjectures 4.1 and 4.2 of \cite{CliftonEtAl}. Moreover, we explicitly compute the mean function on $[0,1]$ of a renewal process with i.i.d. beta($1,m$)-jumps. The expression of this mean function is given in terms of the complex roots of the exponential polynomial $E_m(x) \coloneqq 1 + x/1! + \dots + x^m/m!$, and its derivatives at $0$ are precisely the moments of these roots, as studied in \cite{Zemyan05}.

	The circular convolution identities for Bernoulli polynomials are closely related to the decomposition
	of a real valued random variable $X$ into its integer part $\lfloor X \rfloor \in \Z$ and its fractional 
	part $ X^\circ \in \T := \R/\Z = [0,1)$:
	\begin{equation}
		\label{Xdec}
		X = \lfloor X \rfloor + X^\circ.
	\end{equation}
	If $\gamma_{1}$ is a random variable with standard exponential distribution, then for each positive real $\lambda$ we have the expansion
	\begin{equation} \label{eq:probint}
		\frac{d}{du} \P( (\gamma_1/\lambda)^\circ \le u ) = \frac{\lambda e^{- \lambda u } } { 1 - e^{-\lambda} } = B(u,-\lambda) = \sum_{n \ge 0} b_n(u)  (-\lambda)^n   .
	\end{equation}
	Here the first two equations hold for all real $\lambda \ne 0$ and $u \in [0,1)$, but the final equality holds with a convergent power series only for $0 < |\lambda|  < 2 \pi $. \cref{Sec:convOfSemigroups} presents a generalization of formula \eqref{eq:probint} with the standard exponential variable $\gamma_1$ replaced by the gamma distributed sum $\gamma_r$ of $r$ independent copies of $\gamma_1$, for a positive integer $r$. 
	This provides an elementary probabilistic interpretation and proof of a formula due to Erd\'{e}lyi, Magnus, Oberhettinger, and Tricomi \cite[Section 1.11, page 30]{MR0058756} relating the {\em Hurwitz-Lerch zeta function} (first studied in \cite{Lerch})
	\begin{equation}
		\Phi(z,s,u) = \sum_{m \geq 0} \frac{z^m}{(u + m)^s}
	\end{equation}
	to Bernoulli polynomials. 
	{Moreover, the expansion \eqref{eq:expansion_gamma} in \cref{prop:wrappedGamma} quantifies how the distribution of the fractional part of a $\gamma_{r,\lambda}$ random variable approaches the uniform distribution on the circle in terms of Bernoulli polynomials, where the latter are viewed as signed measures on the circle.}

	\section{Circular convolution of polynomials} \label{Sec:CircularConvOfPoly} 
	\cref{thm:Bernoulliconv} follows easily by induction on $n$ from the characterization \eqref{eq:intchar} of the Bernoulli polynomials, and
	the action of circular convolution by the function
	\begin{equation}
		-b_1(u) = 1/2 - u,
	\end{equation}
	as described by the following lemma.
	\begin{lemma}\label{lem:Thm1Proof}
		For each Riemann integrable function $f$ with domain $[0,1)$, the circular convolution $h = f \circst (-b_1)$ is continuous on $\T$, implying $h(0) = h(1-)$. Moreover, 
		\begin{equation} \label{eq:int0}
			\int_0^1 h(u) du = 0   
		\end{equation}
		and at each $u \in (0,1)$ at which $f$ is continuous, $h$ is differentiable with
		\begin{equation} \label{eq:ddu}
			\frac{d}{du} h(u) = f(u) - \int_0^1 f(v) dv  .
		\end{equation}
		In particular, if $f$ is bounded and continuous on $(0,1)$, then
		$h = f \circst (-b_1)$ is the unique continuous function $h$ on $\T$ 
		subject to \eqref{eq:int0} with derivative \eqref{eq:ddu} at every $u \in (0,1)$.
	\end{lemma}

	\begin{proof}
		According to the definition of circular convolution \eqref{eq:circonv},
		\begin{equation*} \label{eq:fgconv}
			(f \circst g )(u) = \int_0^u f(v) g (u-v) dv + \int_u^1 f(v) g(1 + u - v ) dv .
		\end{equation*}
		In particular, for $g(u) = - b_1(u)$, and  a generic integrable function $f$,
		\begin{align*}
			( f \circst (-b_1) )(u) &= \int_0^u f(v) ( v - u + 1/2) dv  + \int_u^1 f(v) ( v - u - 1/2) dv  \\
			&= \frac{1}{2} \left[ \int_0^u f(v) d v - \int_u ^1 f(v) dv \right] - u \int_0^1 f(v) dv  + \int_0^1 v f(v) dv.
		\end{align*}
		Differentiate this identity with respect to $u$, 
		to see that $h := f \circst (-b_1)$ has 
		the derivative displayed in \eqref{eq:ddu}
		at every $u \in (0,1)$ at which $f$ is continuous, by the fundamental theorem of calculus.
		Also, this identity shows $h$ is continuous on $(0,1)$ with $h(0) = h(0+) = h(1-)$, hence $h$ is continous with respect to the topology of the circle $\T$.
		This $h$ has integral $0$ by associativity of circular convolution: $h \circst 1 = f \circst (-b_1) \circst 1 = f \circst 0 = 0$.
		Assuming further that  $f$ is bounded and continuous on $(0,1)$, the uniqueness of $h$ is obvious.
	\end{proof}
	
	The reformulation of \cref{thm:Bernoulliconv} in \cref{cor:Bernoulliconv} displays how simple it is to convolve  Bernoulli polynomials on the circle.  On the other hand, convolving monomials is less pleasant, as the following calculations show.
	\begin{lemma}\label{lem:convMonomRec}
		For real parameters $n >0$ and $m > - 1$,
		\begin{equation} \label{eq:monoconv}
			{x^m} \circst  {x^n} = {x^n} \circst  {x^m} = \frac{n}{m+1}  {x^{n-1}} \circst  {x^{m+1}} + \frac{ {x^n} -  {x^{m+1}} }{m+1}.
		\end{equation}
	\end{lemma}
	\begin{proof}
		Integrate by parts to obtain
		\begin{align*}
			{x^n} \circst  {x^m} &= \int_{0}^{x} u^{n} (x - u)^m du +  \int_{x}^{1} u^{n} (1 + x - u)^m du  \\
			&= \frac{n}{m+1}\int_{0}^{x} u^{n-1} (x - u)^{m+1} du + \frac{n}{m+1} \int_{x}^{1} u^{n-1} (1 + x - u)^m du + \frac{x^n - x^{m+1}}{m+1} 
		\end{align*}
		and hence \eqref{eq:monoconv}.
	\end{proof}
	
	\begin{proposition}[Convolving monomials]\label{prop:convMonomRec}
		For each positive integer $n$
		\begin{equation}
			\label{eq:ConvMonomials2} {1} \circst   {x^n} = {x^n} \circst  {1}  =  \frac{1}{n+1},
		\end{equation}
		and for all positive integers $m$ and $n$
		\begin{equation} 
			\label{eq:ConvMonomials} {x^m} \circst  {x^n} = {x^n} \circst  {x^m} = \frac{n! \ m!}{(n+m+1)!} + \sum_{k = 0}^{n-1} \frac{n!}{ (n-k)! (m+1)_{k+1}} ( {x}^{n-k} -  {x}^{m+k+1})
		\end{equation}
		and
		with the Pochhammer notation $(m+1)_{k+1} \coloneqq (m+1) \dots (m+k+1)$. In particular
		\[
		{x} \circst  {x^n} = \frac{ {x} -  {x}^{n+1} }{n+1} + \frac{1}{(n+1)(n+2)}.
		\]
	\end{proposition}
	\begin{proof}
		By induction, using \cref{lem:convMonomRec}.
	\end{proof}
	
	\begin{remark} \label{rem:}
		\begin{enumerate}
			\item By inspection of \eqref{eq:ConvMonomials}  the polynomial $\left({x^n} \circst  {x^m}  - \frac{n! \ m!}{(n+m+1)!}\right) / x$ is an anti-reciprocal polynomial with rational coefficients.
			
			\medskip
			
			\item  \cref{thm:Bernoulliconv} can be proved by inductive application of \cref{prop:convMonomRec} to the expansion of the Bernoulli polynomials $B_n(x)$ in the monomial basis. This argument is unnecessarily complicated, but boils down to two following identities for the Bernoulli numbers $B_n:= B_n(0)$ for $n \geq 1$:
			
			\medskip
			
			\begin{align}
				\label{eq:b1id}            B_{n}        &= \frac{-1}{n + 1} \sum_{k = 0}^{n-1} \binom{n+1}{k} B_{k}  \\ 
				\label{eq:b2id} \frac{B_{n+1}}{(n+1)!} &= - \sum_{k = 0}^{n} \frac{1}{(k+2)!} \frac{B_{n-k}}{(n-k)!}.
			\end{align}
			The identity \eqref{eq:b1id} is a commonly used recursion for the Bernoulli numbers. We do not know any reference for \eqref{eq:b2id}, but this can be checked by manipulation of Euler's generating function \eqref{eq:Euler_Gen_Func}. We refer the reader to \cref{appendix:combinatorial proof} for more details.
			
			\medskip
			
			\item Using the hypergeometric function $F \coloneqq \prescript{}{2}{F}_1$, it follows from \cref{eq:ConvMonomials} that:
			$$ x^{n} \circst x^m = \frac{n! m!}{(m+n+1)!} x^{m+n+1} +  \frac{x^{n}}{m+1} F\left(1,-n;m+2; \frac{-1}{x}\right) - \frac{x^{m+1}}{m+1} F(1,-n;m+2;-x).$$
		\end{enumerate}
	\end{remark}

	\section{Probabilistic interpretation} \label{Sec:Prob_Interpretations}
	
	For positive real numbers $a,b > 0$, recall that the beta$(a,b)$ probability distribution, has density
	\[
	\frac{\Gamma(a)\Gamma(b)}{\Gamma(a+b)} \ x^{a-1} (1-x)^{b-1}, \quad ( 0 \leq x  \leq 1)
	\]
	with respect the the Lebesgue measure on $\R$, where $\Gamma$ denotes Euler's gamma function \cite{ArtinGamma} :
	\[
	\Gamma(x) = \int_{0}^{\infty} t^{x-1} e^{-t} dt, \quad \text{for } x > 0.
	\]
	The following corollary offers a probabilistic interpretations of \cref{thm:Bernoulliconv} in terms of the fractional part of a sum of $n$ i.i.d beta$(1,2)$-distributed random variables on the circle.
	
	\begin{corollary}\label{cor:sumBetasModZ}
		The probability density of the sum of $n$ independent beta$(1,2)$ random variables in the circle $\T = \R / \Z$ is
		\[
		(1 - 2 b_1)^{\circst n}(u) = 1 - 2^n b_n(u), \quad \text{for } u \in \T = [0,1).
		\]
	\end{corollary}
	\begin{proof}
		Note that $b_0(u) = 1$ and the density of a beta$(1,2)$ random variable is $2(1 - u) = 1 - 2b_1(u)$ for $0 < u < 1$. So the result follows by induction from \cref{cor:Bernoulliconv}.
	\end{proof}
	
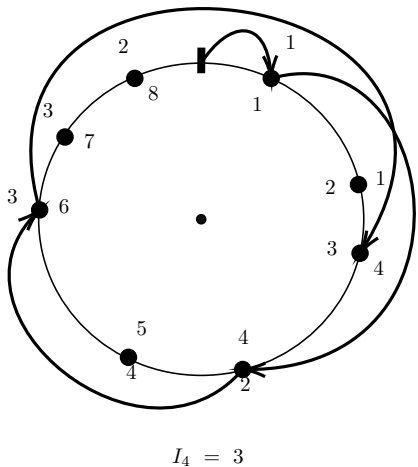
\begin{figure}[ht]
		\begin{center}
			
			\scalebox{0.8}{
				\tikzset{every picture/.style={line width=0.75pt}}      
				\begin{tikzpicture}[x=0.75pt,y=0.75pt,yscale=-1,xscale=1]
					\draw   (178.38,166.55) .. controls (178.38,112.14) and (224.26,68.03) .. (280.86,68.03) .. controls (337.46,68.03) and (383.34,112.14) .. (383.34,166.55) .. controls (383.34,220.96) and (337.46,265.07) .. (280.86,265.07) .. controls (224.26,265.07) and (178.38,220.96) .. (178.38,166.55) -- cycle ;
					\draw [line width=3.75]    (280.86,74.42) -- (280.97,58.25) ;
					\draw  [fill={rgb, 255:red, 15; green, 15; blue, 15 }  ,fill opacity=1 ] (278.01,166.55) .. controls (278.01,164.98) and (279.28,163.7) .. (280.86,163.7) .. controls (282.43,163.7) and (283.71,164.98) .. (283.71,166.55) .. controls (283.71,168.12) and (282.43,169.4) .. (280.86,169.4) .. controls (279.28,169.4) and (278.01,168.12) .. (278.01,166.55) -- cycle ;
					\draw [line width=1.5]    (280.86,68.03) .. controls (297.73,41.68) and (321.76,37.85) .. (324.81,74.75) ;
					\draw [shift={(325,77.67)}, rotate = 267.14] [color={rgb, 255:red, 0; green, 0; blue, 0 }  ][line width=1.5]    (14.21,-4.28) .. controls (9.04,-1.82) and (4.3,-0.39) .. (0,0) .. controls (4.3,0.39) and (9.04,1.82) .. (14.21,4.28)   ;
					\draw [line width=1.5]    (325,77.67) .. controls (431.63,47.15) and (463.84,261.84) .. (309.34,261.36) ;
					\draw [shift={(307,261.33)}, rotate = 0.97] [color={rgb, 255:red, 0; green, 0; blue, 0 }  ][line width=1.5]    (14.21,-4.28) .. controls (9.04,-1.82) and (4.3,-0.39) .. (0,0) .. controls (4.3,0.39) and (9.04,1.82) .. (14.21,4.28)   ;
					\draw [line width=1.5]    (307,261.33) .. controls (244.63,333.93) and (112.67,226.85) .. (176.99,162.6) ;
					\draw [shift={(179,160.67)}, rotate = 137.15] [color={rgb, 255:red, 0; green, 0; blue, 0 }  ][line width=1.5]    (14.21,-4.28) .. controls (9.04,-1.82) and (4.3,-0.39) .. (0,0) .. controls (4.3,0.39) and (9.04,1.82) .. (14.21,4.28)   ;
					\draw [line width=1.5]    (179,160.67) .. controls (124.17,-34) and (496.17,11) .. (381,188) ;
					\draw [shift={(381,188)}, rotate = 303.05] [color={rgb, 255:red, 0; green, 0; blue, 0 }  ][line width=1.5]    (14.21,-4.28) .. controls (9.04,-1.82) and (4.3,-0.39) .. (0,0) .. controls (4.3,0.39) and (9.04,1.82) .. (14.21,4.28)   ;
					\draw  [fill={rgb, 255:red, 0; green, 0; blue, 0 }  ,fill opacity=1 ] (174,160.67) .. controls (174,157.91) and (176.24,155.67) .. (179,155.67) .. controls (181.76,155.67) and (184,157.91) .. (184,160.67) .. controls (184,163.43) and (181.76,165.67) .. (179,165.67) .. controls (176.24,165.67) and (174,163.43) .. (174,160.67) -- cycle ;
					\draw  [fill={rgb, 255:red, 0; green, 0; blue, 0 }  ,fill opacity=1 ] (190,114.67) .. controls (190,111.91) and (192.24,109.67) .. (195,109.67) .. controls (197.76,109.67) and (200,111.91) .. (200,114.67) .. controls (200,117.43) and (197.76,119.67) .. (195,119.67) .. controls (192.24,119.67) and (190,117.43) .. (190,114.67) -- cycle ;
					\draw  [fill={rgb, 255:red, 0; green, 0; blue, 0 }  ,fill opacity=1 ] (230,253.67) .. controls (230,250.91) and (232.24,248.67) .. (235,248.67) .. controls (237.76,248.67) and (240,250.91) .. (240,253.67) .. controls (240,256.43) and (237.76,258.67) .. (235,258.67) .. controls (232.24,258.67) and (230,256.43) .. (230,253.67) -- cycle ;
					\draw  [fill={rgb, 255:red, 0; green, 0; blue, 0 }  ,fill opacity=1 ] (302,261.33) .. controls (302,258.57) and (304.24,256.33) .. (307,256.33) .. controls (309.76,256.33) and (312,258.57) .. (312,261.33) .. controls (312,264.09) and (309.76,266.33) .. (307,266.33) .. controls (304.24,266.33) and (302,264.09) .. (302,261.33) -- cycle ;
					\draw  [fill={rgb, 255:red, 0; green, 0; blue, 0 }  ,fill opacity=1 ] (376,188) .. controls (376,185.24) and (378.24,183) .. (381,183) .. controls (383.76,183) and (386,185.24) .. (386,188) .. controls (386,190.76) and (383.76,193) .. (381,193) .. controls (378.24,193) and (376,190.76) .. (376,188) -- cycle ;
					\draw  [fill={rgb, 255:red, 0; green, 0; blue, 0 }  ,fill opacity=1 ] (320,77.67) .. controls (320,74.91) and (322.24,72.67) .. (325,72.67) .. controls (327.76,72.67) and (330,74.91) .. (330,77.67) .. controls (330,80.43) and (327.76,82.67) .. (325,82.67) .. controls (322.24,82.67) and (320,80.43) .. (320,77.67) -- cycle ;
					\draw  [fill={rgb, 255:red, 0; green, 0; blue, 0 }  ,fill opacity=1 ] (234,77.67) .. controls (234,74.91) and (236.24,72.67) .. (239,72.67) .. controls (241.76,72.67) and (244,74.91) .. (244,77.67) .. controls (244,80.43) and (241.76,82.67) .. (239,82.67) .. controls (236.24,82.67) and (234,80.43) .. (234,77.67) -- cycle ;
					\draw  [fill={rgb, 255:red, 0; green, 0; blue, 0 }  ,fill opacity=1 ] (375,144.67) .. controls (375,141.91) and (377.24,139.67) .. (380,139.67) .. controls (382.76,139.67) and (385,141.91) .. (385,144.67) .. controls (385,147.43) and (382.76,149.67) .. (380,149.67) .. controls (377.24,149.67) and (375,147.43) .. (375,144.67) -- cycle ;
					\draw (332.58,49.3) node [anchor=north west][inner sep=0.75pt]    {$1$};
					\draw (389.33,134.95) node [anchor=north west][inner sep=0.75pt]    {$1$};
					\draw (232,257.07) node [anchor=north west][inner sep=0.75pt]    {$4$};
					\draw (157.11,145.95) node [anchor=north west][inner sep=0.75pt]    {$3$};
					\draw (179.45,91.76) node [anchor=north west][inner sep=0.75pt]    {$3$};
					\draw (304,264.73) node [anchor=north west][inner sep=0.75pt]    {$2$};
					\draw (226.97,51.65) node [anchor=north west][inner sep=0.75pt]    {$2$};
					\draw (388,191.4) node [anchor=north west][inner sep=0.75pt]    {$4$};
					\draw (261,309.73) node [anchor=north west][inner sep=0.75pt]    {$I_{4} \ =\ 3$};
					\draw (311.58,88.3) node [anchor=north west][inner sep=0.75pt]    {$1$};
					\draw (357.58,140.3) node [anchor=north west][inner sep=0.75pt]    {$2$};
					\draw (358.58,179.3) node [anchor=north west][inner sep=0.75pt]    {$3$};
					\draw (302.58,235.3) node [anchor=north west][inner sep=0.75pt]    {$4$};
					\draw (238.58,229.3) node [anchor=north west][inner sep=0.75pt]    {$5$};
					\draw (189.58,152.3) node [anchor=north west][inner sep=0.75pt]    {$6$};
					\draw (205.58,111.3) node [anchor=north west][inner sep=0.75pt]    {$7$};
					\draw (246,81.07) node [anchor=north west][inner sep=0.75pt]    {$8$};
				\end{tikzpicture}
			}
		\end{center}
		\caption{The clock is a circle of circumference $1$. 
			Inside the circle the numbers $1,2, \ldots, 8$ index the order statistics of $8$ uniformly distributed random points on the circle. 
			The corresponding numbers outside the circle are a random assignment of labels from the multiset of four pairs $1^2 2^2 3^2 4^2$. The four successive arrows delimit segments of $\T \equiv [0,1)$ whose lengths $X_1,X_2,X_3,X_4$ are independent beta$(1,2)$ random variables,
			while $(I_1,I_2,I_3,I_4)$ is the sequence of indices inside circle, at the end points of these four arrows. In this example, $(I_1,I_2,I_3,I_4) = (1,4,6,3)$, and the number of turns around the circle is $D_4 = 1$.}
		\label{fig:BernoulliClock}
	\end{figure}

	Recall that a beta($1,2$) random variable can be constructed as the minimum of two independent uniform random variables in $[0,1]$. Let $U_1, U_{1}', \dots ,U_n,U_{n}'$ be a sequence of $2n$ i.i.d random random variables with uniform distribution on $\T = [0,1)$. We think of these variables as random positions around a circle of circumference $1$. On the event of probability one that the $U_i$ and $U_i'$ are all distinct, we define the following variables:
	\begin{enumerate}
		\item $U_{1:2n} < U_{2:2n} < \dots < U_{2n:2n}$ the order statistics of the variables $U_1, U_1', \dots, U_n, U_n'$,
		\item $X_1 := \min(U_1, U_1') $
		\item for $2 \le k \le n$, the variable $X_{k}$ is the spacing around the circle from $X_{k-1}$ to whichever of $U_{k}, U_{k}'$ is encountered first moving clockwise around $\T$ from $X_{k-1}$,
		\item $I_k$ is the random index in $\{1, \dots, 2n\}$ such that $X_k = U_{I_k:2n}$.
		\item $D_n \in \{0, \dots, n-1 \}$ is the random number of full rotations around $\T$ to find $X_n$. This is also the number of descents in the sequence $(I_1,I_2, \dots, I_n)$; that is 
		\begin{equation}
			\label{eq:dndef}
			D_n = \sum_{i = 1}^{n-1} 1[I_{i} > I_{i+1}].
		\end{equation}
	\end{enumerate}
	We refer to this construction as the \emph{Bernoulli clock}. \cref{fig:BernoulliClock} depicts an instance of the Bernoulli clock for $n=4$.

	\begin{proposition} \label{prp:ind}
		With the above notation, the following hold
		\begin{enumerate}
			\item The random spacings $X_1, X_2, \dots , X_n$ (defined by the Bernoulli clock above) are i.i.d beta($1,2$) random variables.
			\item The random sequence of indices $(I_1, I_2, \dots, I_n)$ is independent of the sequence of order statistics $(U_{1:2n} , \dots , U_{2n:2n})$.
		\end{enumerate}
	\end{proposition}
	\begin{proof}
		This is a corollary of \cref{prop:BernoulliClockGenCase}. See \cref{Sec:Generalized_Bernoulli_clock} where the general case is discussed.
	\end{proof}
	
	\subsection{Expanding Bernoulli polynomials in the Bernstein basis}
	
	It is well known that, for $1 \leq k \leq 2n $ the distribution of $U_{k:2n}$ is beta($k, 2n+1-k$), whose
	probability density relative to Lebesgue measure at $u \in [0,1)$ is the normalized Bernstein polynomial of degree $2 n - 1$
	\[
	f_{k:2n}(u) := \frac{(2n)!}{(k-1)!(2n - k)!} u^{k - 1}(1-u)^{2n - k}  
	\]
	\begin{proposition}\label{prop:ExpansionInBernstein}
		For each positive integer $n$, the sum $S_n$ of $n$ independent beta$(1,2)$ variables has fractional part $S_n^\circ$
		whose probability density on $(0,1)$ is given by the formulas
		\begin{equation} \label{eq:expandBetainBernstein}
			f_{S_n}^\circ(u) = 1 - 2^n b_n(u) = \sum_{k = 1}^{2n} p_{k:2n} \ f_{k:2n}(u),  \quad \text{for } u \in (0,1).
		\end{equation}
		where $(p_{1:2n}, \dots, p_{2n:2n})$ is the probability distribution of the random index $I_n$ in the Bernoulli clock construction:
		\[
		p_{k:2n} = \P(I_n = k), \quad \text{for } 1 \leq k \leq 2n.
		\]
	\end{proposition}
	\begin{proof}
		The first formula for the density of $S_n^\circ$ is read from \cref{cor:sumBetasModZ}.
		Proposition \ref{prp:ind} represents $S_n^\circ = U_{I_n: 2n}$ where the index $I_n$ is independent of the sequence of
		order statistics $(U_{k:2n}, 1 \le k \le 2 n)$, hence the second formula for the same probability density on $(0,1)$.
	\end{proof}
	
	\begin{corollary}
		The factorially normalized Bernoulli polynomial of degree $n$ admits the expansion in Bernstein polynomials of degree $2 n - 1$
		\begin{equation} \label{eq:expandingBernoulliInBernstein}
			b_{n}(u) = \frac{1}{2^n} \sum_{k = 1}^{2n} \delta_{k: 2n } \ f_{k:2n}(u)
		\end{equation}
		where $\delta_{k:2n}$ is the difference at $k$ between the uniform probability distribution on $\{1, \dots, 2n \}$ and the distribution of $I_n$.
		\begin{equation}
			\label{eq:deldef}
			\delta_{k:2n} = \frac{1}{2n} - p_{k:2n}  \quad \text{for } 1 \leq k \leq 2n. 
		\end{equation}
	\end{corollary}
	\begin{proof}
		Formula \eqref{eq:expandingBernoulliInBernstein} is obtained from \eqref{eq:expandBetainBernstein}i,
		in the first instance as an identity of continuous functions of $u \in (0,1)$, then as an identity of polynomials in $u$, 
		by virtue of the binomial expansion 
		\[
		\sum_{k=1}^{2n} \frac{1}{2n} f_{k:2n}(u) = 1. \qedhere
		\]
	\end{proof}
	
	\begin{remark}
		Since $b_n(1 - u) = (-1)^{n} b_{n}(u)$ and $f_{k:2n}(1-u) = f_{2 n + 1 - k: 2n} (u)$, the identity \eqref{eq:expandingBernoulliInBernstein} implies that the
		difference between the distribution of $I_n$ and the uniform distribution on $\{1, \ldots, 2 n\}$ has the symmetry
		\begin{equation} \label{eq:deltaSymm}
			\delta_{2n + 1 - k : 2n} = (-1)^n \delta_{k:2n}   \quad \text{for } 1 \leq k \leq 2n.
		\end{equation}
	\end{remark}
	
	\begin{conjecture} We conjecture that the discrete sequence $(\delta_{1:2n}, \dots, \delta_{2n:2n})$ approximates the Bernoulli polynomials $b_n$ (hence also the shifted cosine functions, see \eqref{eq:SineFunctionApprox}) as $n$ becomes large, more precisely:
		$$ \sup\limits_{1 \leq k \leq 2n} \left| 2n \pi^n \delta_{k:2n}- (2\pi)^n b_n\left( \frac{k-1}{2n-1}\right) \right|  \xrightarrow[]{} 0 \quad \text{as } n \to \infty.$$
	\end{conjecture}
	
\begin{figure}[ht]
		\begin{center}
			\centering
			\includegraphics[scale=0.30]{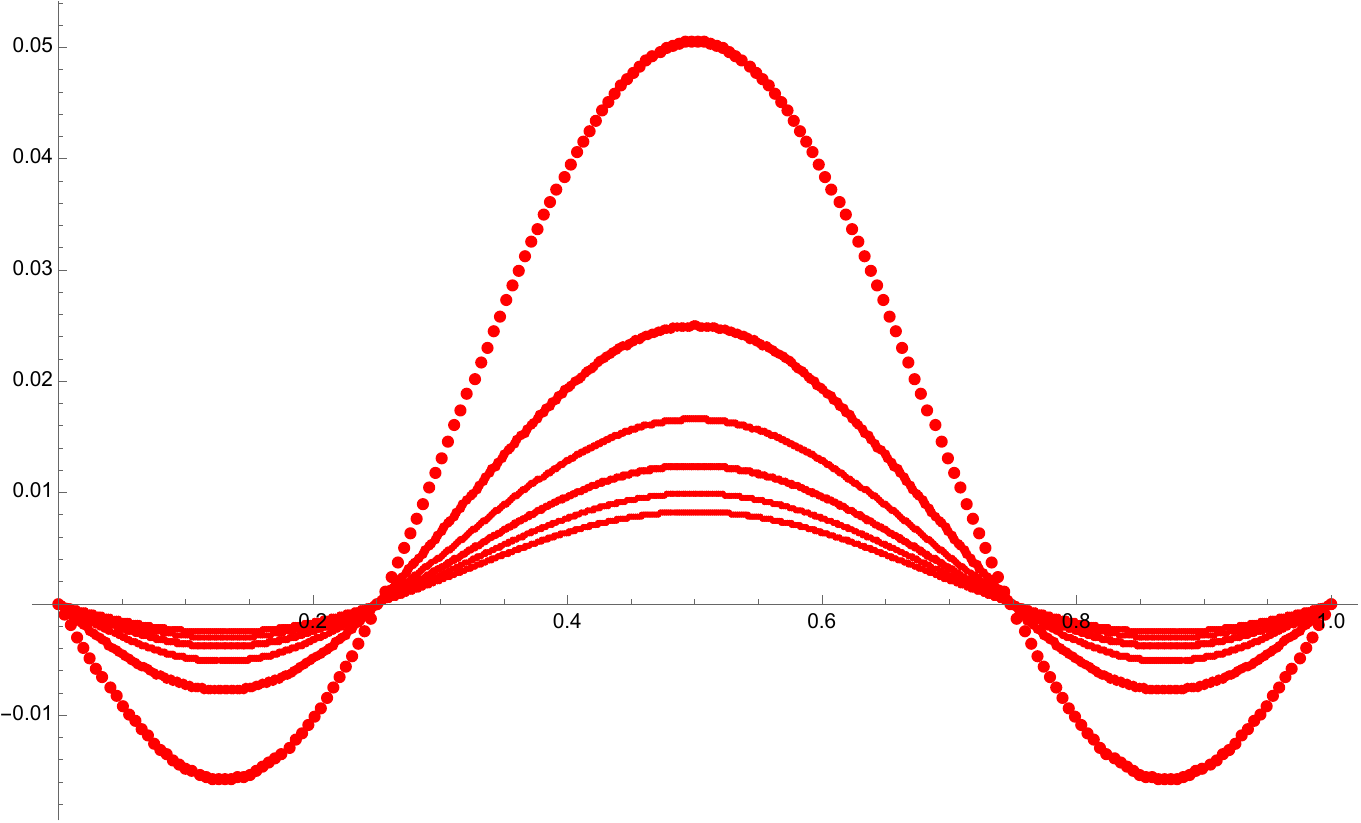}
		\end{center}
		\caption{Plots of $2n \pi^n \delta_{k:2n} - (2\pi)^n b_n \left( \frac{k-1}{2n -1 }\right)$ for $n = 100, 200, 300, 400, 500, 600$.}
		\label{fig:diff_plots}
	\end{figure}

	\cref{fig:diff_plots} does suggest that the difference $2n \pi^n \delta_n(k) - (2\pi)^n b_n\left( \frac{k-1}{2n-1}\right)$ gets smaller uniformly in $1 \leq k \leq 2n$ as $n$ grows, geometrically but rather slowly, like $C \rho^n$ for a constant $C >0$ and $\rho \approx 2^{- 1/100}$.
	\begin{figure}[ht] \label{fig:delta_plots}
		\centering
		\begin{subfigure}[b]{0.35\textwidth}
			\centering
			\includegraphics[scale=0.35]{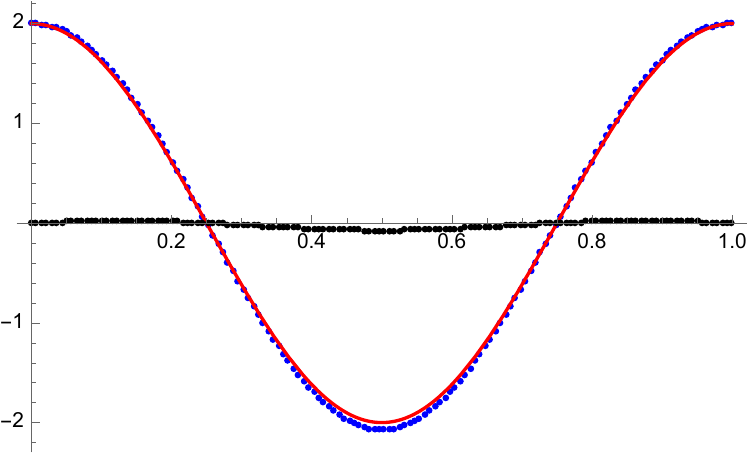}
			\caption[$n=70$]%
			{{\small $n=70$ }}    
		\end{subfigure}
		\hfill
		\begin{subfigure}[b]{0.35\textwidth}  
			\centering 
			\includegraphics[scale=0.35]{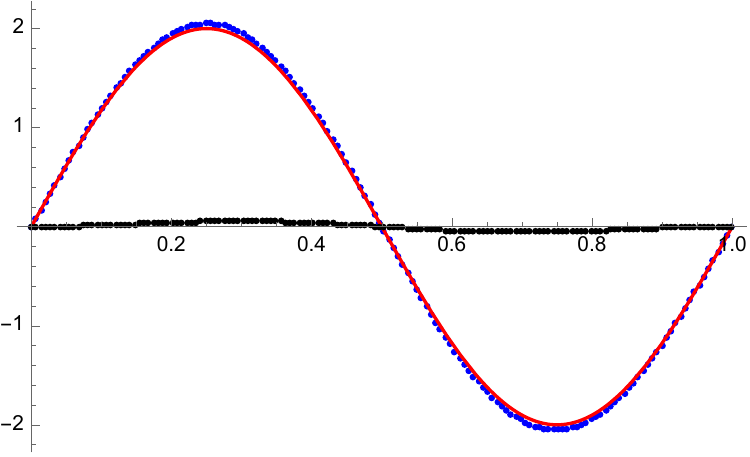}
			\caption[$n=75$]%
			{{\small $n=75$}}    
		\end{subfigure}
		\vskip\baselineskip
		\begin{subfigure}[b]{0.35\textwidth}   
			\centering 
			\includegraphics[scale=0.35]{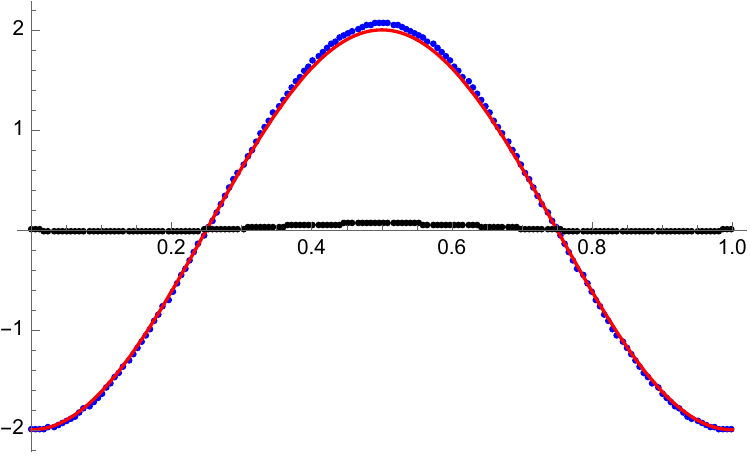}
			\caption[$n=80$]
			{{\small $n=80$}}    
			\label{fig:mean and std of net34}
		\end{subfigure}
		\hfill
		\begin{subfigure}[b]{0.35\textwidth}   
			\centering 
			\includegraphics[scale=0.35]{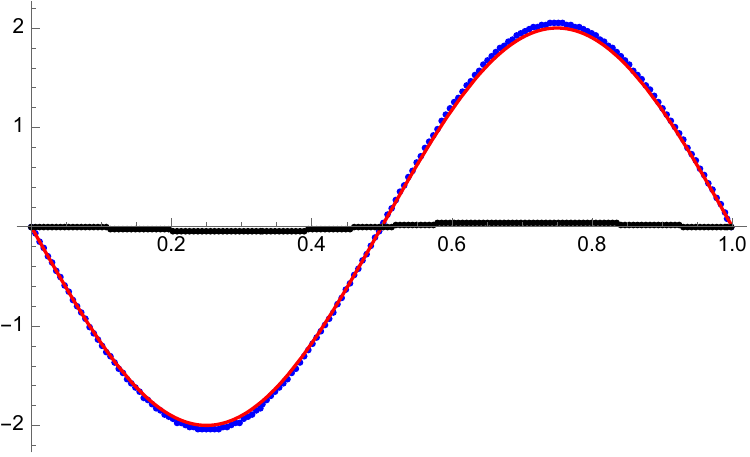}
			\caption[$n=85$]%
			{{\small $n=85$}}    
		\end{subfigure}
		\caption{Plots of $2n \pi^n \delta_n$ (dotted curve in blue), $(2\pi)^n b_n(x)$ (curve in red) and their difference (dotted curve in black) for $n = 70, 75, 80, 85$.}
	\end{figure}
	
	From \eqref{eq:expandBetainBernstein} we see that we can expand the polynomial density $1 -2^n b_n(u)$ in the Bernstein basis of degree $2n - 1$ with positive coefficients. A similar expansion can obviously be achieved using Bernstein polynomials of degree $n$, with coefficients which must add to $1$. These coefficients are easily calculated for modest values of $n$ (see \eqref{eq:bernoulliExpansionBernsteine}) which suggests the following
	
	\begin{conjecture}
		For each positive integer $n$, the polynomial probability density $1 -2^n b_n(u)$ on $[0,1)$ can be expanded in the Bernstein basis of degree $n$ with positive coefficients.
	\end{conjecture}
	\begin{question}
		More generally, what can be said about the greatest multiplier $c_n$ such that the polynomial $1 - c_n b_n(x)$ is a linear combination of degree $n$ Bernstein polynomials with non-negative coefficients?
	\end{question}
	
	\subsection{The distributions of \texorpdfstring{$I_n$}{} and \texorpdfstring{$D_n$}{}}
	
	\begin{proposition}
		The distribution of $I_n$ in the Bernoulli clock construction is given by 
		\begin{equation}
			\label{eq:p_nExplicit}
			\P( I_n = k ) = \frac{1}{2n} - \delta_{k:2n} \quad \text{for } 1 \leq k \leq 2n \mbox{ with }
		\end{equation}
		\begin{equation} \label{eq:delta_nExplicit}
			\delta_{k:2n} = \frac{2^{n-1}}{n \ n!}\sum_{i = 0}^{n}  \frac{\binom{k-1}{i} \binom{n}{i} }{ \binom{2n - 1}{i}}  B_{n-i} ,  \quad \text{ for } 1 \leq k \leq 2n.
		\end{equation}
	\end{proposition}
	\begin{proof}
		For each positive integer $N$, in the Bernstein basis $(f_{j:N})_{1 \leq j \leq N}$ of polynomials of degree at most $N-1$, it is well known that the monomial $x^i$ can be expressed as
		\[
		x^i = \frac{1}{N \binom{N - 1}{i}} \sum_{j = i+1}^{N} \binom{j-1}{i} f_{j:N}(x) \quad \text{for }  0 \leq i < N,
		\]
		see \cite[Table 2.1]{Riordan} for a reference.
		Plugging this expansion into \eqref{eq:bernexp} yields the expansion of $b_n(x)$ in the Bernstein basis of degree $N-1$ for every $N > n$
		\begin{equation}\label{eq:bernoulliExpansionBernsteine}
			b_n(x) = \sum_{j = 1}^N  \left(\sum_{i = 0}^{n}  \frac{\binom{j-1}{i} \binom{n}{i}}{ n! N \binom{N - 1}{i}}  B_{n-i} \right)  f_{j:N}(x) \qquad ( 0 \le n < N).
		\end{equation}
		
		In particular, for $N = 2 n$ comparison of this formula with \eqref{eq:expandingBernoulliInBernstein} yields 
		\eqref{eq:delta_nExplicit} and hence \eqref{eq:p_nExplicit}
	\end{proof}
	
	\begin{remark}
		The error $\delta_{k:2n}$ is polynomial in $k$ and the symmetry $\delta_{2n + 1 - j : 2n} = (-1)^n \delta_{j : 2n}  $ is not obvious from \eqref{eq:delta_nExplicit}.
	\end{remark}
	
	Let us now derive the distribution of $D_n$ explicitly. From the Bernoulli clock scheme, we can construct the random variable $D_n$ as follows. Let $X_1, \dots, X_n$ be a sequence of i.i.d random variables and $S_n := X_1 + \dots + X_n$ their sum in $\R$ (not in the circle $\T$), then
	\[
	D_n = \lfloor S_n \rfloor.
	\]
	
	\begin{theorem}
		The distribution function of $S_n$ is given by
		\[
		\P(S_n \leq x) = 2^n \sum_{k =0}^{n}\sum_{ j = 0}^{n-k}  \binom{n}{k} \binom{n-k}{j} (-1)^{n-k-j} \frac{(x-k)_+^{2n - j}}{(2n-j)!}, \quad \quad \text{for } x \geq 0,
		\]
		where $x_{+}$ denotes $\max(x,0)$ for $x \in \R$.
	\end{theorem}
	\begin{proof}
		Let $\varphi$ be the Laplace transform of the $X_i$'s i.e.
		\[
		\varphi_X(\theta) \coloneqq \E[e^{- \theta X_1}] = \int_{0}^{+ \infty} \theta e^{-\theta t} \P(X_1 \leq t) dt, \quad \text{for } \theta > 0.
		\]
		We compute $\varphi_{X}$ and we obtain
		\[
		\varphi_X(\theta) = \frac{2}{\theta^2} \left( e^{- \theta} + (\theta - 1) \right), \quad \text{for } \theta > 0.
		\]
		So for $n \geq 1$, the Laplace transform of $S_n$ is then given by
		\begin{equation}\label{eq:LaplaceTransformS_n}
			\varphi_{S_n}(\theta) = \left(\varphi_X(\theta)\right)^n = 2^n \sum_{k =0}^{n}\sum_{ j = 0}^{n-k}  \binom{n}{k} \binom{n-k}{j} (-1)^{n-k-j} \frac{ e^{- k \theta} }{\theta^{2n - j}}.
		\end{equation}
		The transform $\varphi_{S_n}$ can be inverted term by term using the following identity
		\begin{equation}\label{eq:masterLaplaceTransform}
			\int_{0}^{+\infty} \theta e^{-\theta t} \frac{(t - k)_{+}^{n}}{n!} \ dt = \frac{e^{- k\theta }}{\theta^n}, \quad \text{for } k \geq 0, \ \theta > 0 \text{ and } n \geq 0.
		\end{equation}
		We then obtain the cdf of $S_n$ as follows:
		\begin{equation}\label{eq:cdfOfS_n}
			\P(S_n \leq x) = 2^n \sum_{k =0}^{n}\sum_{ j = 0}^{n-k}  \binom{n}{k} \binom{n-k}{j} (-1)^{n-k-j} \frac{(x-k)_+^{2n - j}}{(2n-j)!}, \quad \text{for } x \geq 0. \qedhere
		\end{equation}
	\end{proof}

	\begin{remark}
		\eqref{eq:masterLaplaceTransform} was known to Lagrange in the 1700s and it appears in \cite[Lemme III and Corollaire I]{lagrangeMemoire} where he said the final words on inverting Laplace transforms of the form (\ref{eq:LaplaceTransformS_n}):
		\begin{quote}
			\textit{"... mais comme cette intégration est facile par les methodes connues, nous n'entrerons pas dans un plus grand detail là-dessus; et nous terminerons même ici nos recherches, par lesquelles on doit voir qu'il ne reste plus de difficulté dans la solution des questions qu'on peut proposer à ce sujet."}
		\end{quote}
	\end{remark}
	
	Since $S_n$ has a density, we can deduce that
	\[
	\P(D_n = k) = \P(S_n \leq k + 1) - \P(S_n \leq k), \quad \text{for } 0 \leq k \leq n - 1.
	\]
	Combined with \eqref{eq:cdfOfS_n} this gives the distribution of $D_n$ explicitly. The following table gives the values of the number of permutations of the multiset $1^2 \dots n^2$ for which $D_n = d$, which we denote by $\#(n; +, d)$, for small values of $n$.
	
	\begin{table}[H]
		\centering
		\begin{tabular}{|l|c|c|c|c|c|c|c|}
			\hline
			\diagbox{$n$}{$d$} & $0$ & $1$ & $2$ & $3$  & $4$ & $5$ \\
			\hline
			$1$	& $1$  &   &   &  &   &  \\
			\hline
			$2$	& $1$ & $ 1 $ &    &   &   &   \\
			\hline
			$3$	& $47$ & $42$ & $1$ &   &  &   \\
			\hline
			$4$	& $641$ & $1659$ & $219$ &  $1$ &  &   \\
			\hline
			$5$	& $11389$ & $72572$ & $28470$ & $968$ & $1$ &  \\
			\hline
			$6$	& $248749$ & $3610485$ & $3263402$ & $357746$ & $4017$ & $1$ \\
			\hline
		\end{tabular}
		\caption{The table of $\#(n; + , d)$.}
		\label{tab:distOfD_n}
	\end{table}
	
	\begin{remark}
		The sequence $a(n) = \#(n; +, 0) = 2^{-n} (2n)! \ \P(D_n = 0)$, which counts the number of permutations of $1^2 \cdots n^2$ for which $D_n = 0$ (the first column in \cref{tab:distOfD_n}), can be explicitly written using \eqref{eq:cdfOfS_n} as follows
		\begin{equation}\label{eq:a(n)Explicit}
			a(n) =\P(S_n \leq 1) =  \sum_{j = 0}^{n} (-1)^{n-j}  \binom{n}{j}\frac{  (2n)!}{(2n-j)!}.
		\end{equation}
		This integer sequence appears in many other contexts (see OEIS entry \href{https://oeis.org/A006902}{A006902}), among which we mention a few:
		\begin{enumerate}
			\item $a(n)$ is the number of words on $1^2 \cdots n^2$ with longest complete increasing sub-sequence of length $n$. We shall detail this in \cref{Sec:Generalized_Bernoulli_clock}.
			
			\item $a(n) = n! \ Z(\mathfrak{S}_n; n, n-1, \dots , 1)$ where $Z(\mathfrak{S}_n)$ is the cycle index of the symmetric group of order $n$ (see \cite[Section 1.3]{StanleyVol1}).
			
			\item $a(n) = \mathrm{B}_n \left(n \cdot 0!, \ (n-1) \cdot 1!, \ (n-2)! \cdot 2!, \ \dots, \ 1 \cdot (n-1)! \right)$, where $\mathrm{B}_n(x_1, \dots ,x_n)$ is the $n$-th complete Bell polynomial.
		\end{enumerate}
		
	\end{remark}
	
	\section{Combinatorics of the Bernoulli clock}\label{Sec:BernoulliClock}
	
	There are a number of known constructions of the Bernoulli numbers $B_n$ by permutation enumerations. Entringer \cite{MR205866} showed that Euler's presentations of the Bernoulli numbers, as coefficients in the expansions of hyperbolic and trigonometric functions, lead to explicit formulas for $B_n$ by enumeration of alternating permutations. More recently, Graham and Zang \cite{graham_enumerating_2008} gave a formula for $B_{2n}$ by enumerating a particular subset of the set of $2^{-n}(2n)!$ permutations of the multiset $1^2 \cdots n^2$ of $n$ pairs.
	
	The number of permutations of this multiset, such that for every $i < n$ between each pair of occurrences of $i$ there is exactly one $i+1$, is $(-2)^n ( 1 - 2^{2 n} )B_{2n}$. Here we offer a novel combinatorial expression of the Bernoulli numbers 
	based on a different attribute of permutations of same multiset \eqref{eq:ourmulti},
	which arises from the the probabilistic interpretation in \cref{Sec:Prob_Interpretations}.
	We call the combinatorial construction involved the \emph{the Bernoulli clock}. Fix a positive integer $n \geq 1$ and for a permutation $\tau$ of the multiset \eqref{eq:ourmulti}, 
	
	\begin{itemize}
		\item Let $1 \leq I_1 \leq 2n-1$ be the position of the first $1$; that is $I_1 = \min\{1  \leq k \leq 2n \colon  \tau(k) = 1 \}$.
		\item For $1 \leq k \leq n-1$, denote by $1 \leq I_{k+1} \leq 2n$ the index of the first value $k+1$ following $I_k$ in the cyclic order (circling back to the beginning of necessary).
		\item Let $0 \leq D_n \leq n-1$ be the number of times we circled back to the beginning of the multiset before obtaining the last index $I_n$.
	\end{itemize}
	
	\begin{example}
		The permutation $\tau$ corresponding to \cref{fig:BernoulliClock} is the permutation $\tau = (1,1,4,2,4,3,3,2)$. For this permutation
		$$
		(I_1,I_2,I_3,I_4) = (1,4,6,3) \quad \text{and} \quad D_4 = 1 .
		$$
	\end{example}

	Notice that random index $I_n$ and the number of descents $D_n$ depend only on the relative positions of $U_1, U_1', \dots, U_{n},U_{n}'$ i.e. the permutation of the multiset $1^2 \cdots n^2$. 
	So the distribution of $I_n$ and $D_n$ can be obtained by enumerating permutations. For $n \geq, 1 \leq i \leq 2n$ and $0 \leq d \leq n-1$, let us denote by
	\begin{enumerate}
		\item $\#(n;i,d)$ the number of permutations among the $(2n)! / 2^n$ permutations of the multiset $\{1,1, \dots, n,n\}$ that yield $I_n = i$ and $D_n = d$,
		\item $\#(n; i, +)$ the number of permutations that yield $I_n = i$,
		\item $\#(n; +, d)$ the number of permutations that yield $D_n = d$.
	\end{enumerate}
	For $n = 2$ there are $6$ permutations of $\{1,1,2,2\}$ summarized in the following table
	\begin{table}[H]
		\begin{center}	
			\begin{tabular}{|c|c|c|c|c|c|c|}
				\hline
				Permutations   & $1122$ & $1212$ & $1221$ & $2112$ & $2121$ & $2211$ \\
				\hline
				($I_2, D_2$) & $(3,0)$ & $(2,0)$ & $(2,0)$ & $(4,0)$ & $(3,0)$ & $(1,1)$\\
				\hline
			\end{tabular}
			\caption{Permutations of $\{1,1,2,2\}$ and corresponding values of $(I_2, D_2)$.}
			\label{tab:permutations(I_2,D_2)}
		\end{center}
	\end{table}
	The joint distribution of $I_2, D_2$ is then given by 
	
	\begin{table}[H]
		\begin{center}	
			\begin{tabular}{|l|c|c|c|c|c|}
				\hline
				\diagbox{$I_2$}{$D_2$}		 & $1$ & $2$ & $3$ & $4$  & $\#(2; +, \bullet)$  \\
				\hline
				$0$	& $0$ & $2$ & $2$ & $1$ & 5  \\
				\hline
				$1$	& $1$ & $0$ & $0$ & $0$ &  1 \\
				\hline
				$\#(2; \bullet, +)$  & $1$ & $2$ & $2$ & $1$ &  6 \\
				\hline
			\end{tabular}
			\caption{The table of $\#(2; \bullet, \bullet)$.}
			\label{tab:distributionOf(I_2,D_2)}
		\end{center}
	\end{table}
	
	Similarly for $n=3$ we get 
	\begin{table}[H]
		\begin{center}	
			\begin{tabular}{|l|c|c|c|c|c|c|c|c|}
				\hline
				\diagbox{$I_3$}{$D_3$} & $1$ & $2$ & $3$ & $4$  & $5$ & $6$ &  $\#(3; +, \bullet)$  \\
				\hline
				$0$	& $0$ & $0$ & $6$ & $12$ & $15$ & $14$ & $47$ \\
				\hline
				$1$	& $14$ & $13$ & $8$ & $4$ & $2$ & $1$& $42$\\
				\hline
				$2$	& $1$ & $0$ & $0$ & $0$ & $0$& $0$ &$1$\\
				\hline
				$\#(3; \bullet, +)$  & $15$ & $13$ & $14$ & $16$  & $17$ & $15$& $90$\\
				\hline
			\end{tabular}
			\caption{The table of $\#(3; \bullet, \bullet)$.}
			\label{tab:distributionOf(H_3,D_3)}	
		\end{center}
	\end{table}
	
	The distribution of $(I_n, D_n)$ can be obtained recursively as follows. The key observation is that every permutation of the multi-set $1^2 2^2 \cdots n^2$ is obtained by first choosing a permutation of  $1^2 2^2 \cdots (n-1)^2$, then choosing $2$ places to insert the two values $n,n$. There are $\binom{2(n-1)}{2}$ options for where to insert the two last values. This corresponds to the factorization 
	\[
	(2n)! \ 2^{-n} = (2(n-1))! \ 2^{-n+1} \  \binom{2n}{2}.
	\]
	Moreover, for $x \in \{1, \dots, 2(n-1)\}$ the identity of quadratic polynomials
	\[
	\binom{x+1}{2} + \binom{2n - x}{2} + x(2n-1 - x) = \binom{2n}{2},
	\] 
	translates, for each integer $x \in \{1, \dots, 2(n-1)\}$ and each permutation $\sigma$ of $1^2,\cdots (n-1)^2$, the decomposition of the total number of ways to insert the next two values $n, n$ according to whether:
	\begin{enumerate}
		\item both places are to the left of $x$,
		\item both places are to the right of $x$,
		\item one of those places is to the left of $x$ and the other to the right of $x$.
	\end{enumerate}
	
	Suppose we ran the Bernoulli clock scheme on $2(n-1)$ hours and obtained $(I_{n-1}, D_{n-1})$. Inserting two new values $n, n$, the index $I_{n}$ then depends only on $I_{n-1}$ and the places where the two new values $n$ are inserted relatively to $I_{n-1}$. So, the sequence $(I_1, I_2, \dots )$ is a time-inhomogeneous Markov chain starting from $I_1 = 1$ and a $2(n-1) \times 2n$ transition matrix from $I_{n-1}$ to $I_n$ given by
	\[
	P_{n}(x \to y)  = \P(I_n = y | I_{n-1}  = x) = \frac{Q_n(x,y)}{\binom{2n}{2}}, \qquad  (1 \leq x \leq (2n-1), \ 1 \leq y \leq 2n)
	\]
	where $Q_n(x,y)$ is the number of ways to insert the two new values $n$ in the Bernoulli clock in such a way that the first one of them to the right of $x$ is at place $y$. More explicitly, by elementary counting, we have
	\[
	Q_n(x,y) = 
	\begin{cases} x - y + 1,   \quad \quad \quad  \ \ \text{if }  1 \leq y \leq x  \\
		2n - 1 - x , \quad \quad \quad  \text{if }  y = x+1 \\ 
		2n - y + x,  \quad \quad \quad \text{if }  x+2 \leq y \leq 2n
	\end{cases}
	\]
	So the first few transition matrices are
	
	\[  P_2 = \frac{Q_2}{\binom{4}{2}} =\frac{1}{6} \left(
	\begin{array}{cccc}
		1 & 2 & 2 & 1 \\
		2 & 1 & 1 & 2
	\end{array}
	\right), \quad \quad
	P_3 = \frac{Q_3}{\binom{6}{2}} = \frac{1}{15}  \left(
	\begin{array}{cccccc}
		1 & 4 & 4 & 3 & 2 & 1 \\
		2 & 1 & 3 & 4 & 3 & 2 \\
		3 & 2 & 1 & 2 & 4 & 3 \\
		4 & 3 & 2 & 1 & 1 & 4
	\end{array}
	\right),
	\]
	\[ \text{and} \quad
	P_4 = \frac{Q_4}{\binom{8}{2}} = \frac{1}{28} \left(
	\begin{array}{cccccccc}
		1 & 6 & 6 & 5 & 4 & 3 & 2 & 1 \\
		2 & 1 & 5 & 6 & 5 & 4 & 3 & 2 \\
		3 & 2 & 1 & 4 & 6 & 5 & 4 & 3 \\
		4 & 3 & 2 & 1 & 3 & 6 & 5 & 4 \\
		5 & 4 & 3 & 2 & 1 & 2 & 6 & 5 \\
		6 & 5 & 4 & 3 & 2 & 1 & 1 & 6
	\end{array}
	\right),
	\]
	see \cref{tab:ConstructionOfP_3} for a detailed combinatorial construction of $Q_3$. This discussion is summarized by the following proposition. 
	
	\begin{proposition} 
		For a uniform random permutation of $1^2 \cdots n^2$ the probability distribution of $I_n$, 
		treated as a $1\times 2n $ row vector $p_n = (p_{1:2n}, \dots, p_{2n:2n})$, 
		is determined recursively by the matrix forward equations 
		\begin{equation}
			\label{eq:forwardEquation} p_{n+1}  = p_{n}  \ P_{n+1} \qquad \mbox{ for } n = 1,2, \ldots \mbox{ starting from } p_1 = (1,0).\\
		\end{equation} 
	\end{proposition}

	\noindent So the first few of these distributions of $I_n$ are as follows:
	\begin{align*}
		p_1 &= (1,0),  								    & p_2 &= \frac{1}{6} (1,2,2,1), \\ 
		p_3 &= \frac{1}{90} (15,13,14,16,17,15),    & p_4 &= \frac{1}{2520} (322,322,312,304,304,312,322,322).
	\end{align*}
	As $n$ become bigger, the distribution $p_{n}$ gets closer to the uniform on $\{1, \dots, 2n\}$. The error $\delta_n(k) = 1/(2n) - p_{k:2n}$ is polynomial in $k$ and satisfies the same forward equation as $p_n$ i.e.
	\begin{equation}
		\label{eq:forwardEquation_for_delta} \delta_{n+1} = \delta_{n}  \ P_{n+1} \qquad \mbox{ for } n = 1,2, \ldots \mbox{ starting from } \delta_0 = (1/2,-1/2).
	\end{equation}
	The sequence $\delta_n$ is also closely tied to the polynomial $b_n(x)$ as \eqref{eq:expandingBernoulliInBernstein} shows.
	
	\begin{example}\label{ex:ConstructionOfP_3}
		Let us detail the combinatorics of permutations that yields the matrix $P_3$.                    
		\begin{table}[H]
			\centering
			\begin{tabular}{|c|c|c|c|c|c|}
				\hline
				1 & 2 & 3 & 4 & 5 & 6\\
				\hline
				\multicolumn{6}{l}{}\\
				\hline
				$\bullet^{1,2,3,4}$ &  $\bullet$        & 1                & 2                & 3                & 4\\
				\hline
				$\bullet^{2,3,4}$  & 1                 & $\bullet^{1}$    & 2                & 3                & 4\\
				\hline
				$\bullet^{3,4}$   & 1                 & 2                &  $\bullet^{1,2}$  & 3                & 4\\
				\hline
				$\bullet^{4}$    & 1                 & 2                & 3                & $\bullet^{1,2,3}$  & 4\\
				\hline
				$\bullet$        & 1                 & 2                & 3                & 4                & $\bullet^{1,2,3,4}$ \\
				\hline
				1                & $\bullet^{1,2,3,4}$  & $\bullet$        & 2                & 3                & 4\\
				\hline
				1                & $\bullet^{1,3,4}$   & 2                & $\bullet^{2}$    & 3                & 4\\
				\hline
				1                & $\bullet^{1,4}$    & 2                & 3                & $\bullet^{2,3}$   & 4\\
				\hline
				1                & $\bullet^{1}$     & 2                & 3                & 4                & $\bullet^{2,3,4}$\\
				\hline
				1                & 2                 & $\bullet^{1,2,3,4}$ & $\bullet$        & 3                & 4\\
				\hline
				1                & 2                 & $\bullet^{1,2,4}$  & 3                & $\bullet^{3}$    & 4\\
				\hline
				1                & 2                 & $\bullet^{1,2}$   & 3                & 4                & $\bullet^{3,4}$\\
				\hline
				1                & 2                & 3                 & $\bullet^{1,2,3,4}$ & $\bullet$        & 4\\
				\hline
				1                & 2                & 3                 & $\bullet^{1,2,3}$  & 4                & $\bullet^{4}$\\
				\hline
				1                & 2                & 3                 & 4                & $\bullet^{1,2,3,4}$ & $\bullet$\\
				\hline
				
				\multicolumn{6}{l}{ }\\
				
				\hline
				$ {1}$ & $ {4}$ & $ {4}$ & $ {3}$ & $ {2}$ & $ {1}$\\
				\hline 
				$ {2}$ & $ {1}$ & $ {3}$ & $ {4}$ & $ {3}$ & $ {2}$\\
				\hline 
				$ {3}$ & $ {2}$ & $ {1}$ & $ {2}$ & $ {4}$ & $ {3}$\\
				\hline 
				$ {4}$ & $ {3}$ & $ {2}$ & $ {1}$ & $ {1}$ & $ {4}$ \\
				\hline
			\end{tabular}                            
			\caption{Combinatorial construction of $Q_3$: The top $1\times 6$ row displays the column index of places in rows of the main $15 \times 6$ table below it. The $15$ rows of the main table list all $\binom{6}{2} = 15$ pairs of places, 
				represented as two dots $\bullet$, in which two new values $3,3$  can be inserted relative to $4$ possible places
				of $I_2 \in \{1, 2, 3, 4\}$.  The exponents of each dot $\bullet$ are the values of $I_2$ leading to $I_3$ being the column index of that dot in $\{1,2, 3,4,5, 6 \}$. For example in the second row, representing insertions of the new value $3$ in places $1$ and $3$ of $6$ places, the dot $\bullet^{2,3,4}$ in place $1$ is the place $I_3$ found by the Bernoulli clock algorithm if $I_2 \in \{2,3,4\}$. The matrix $Q_3$ is the $4 \times 6$ matrix below the main table.
				The entry $Q_3(i,j)$ in row $i$ and column $j$ of $Q_3$ is the number of times $i$ appears in the exponent of a dot $\bullet$ in the $j$-th column of the main table.}
			\label{tab:ConstructionOfP_3}
		\end{table}
		
		Notice that the matrices $Q_n$ have the remarkable symmetry
		\begin{equation}\label{eq:SymmOfQ}
			2n - 1 - Q_{n}(i,j) = \widetilde{Q}_{n}(i, j), \quad (1 \leq i \leq 2n, \ 1 \leq j \leq 2n+2),
		\end{equation}
		with $\widetilde{Q}_{n}(i, j) \coloneqq Q_{n}(2n - 1 - i, 2n + 1  -j)$ i.e. the matrix $\widetilde{Q}_n$ is the matrix $Q_n$ with entries in reverse order in both axis.
	\end{example}
	
	\vspace{2mm}
	
	\begin{remark}~\\
		\begin{enumerate}
			\item It is interesting to note that, from \eqref{eq:forwardEquation}, it is not clear what the Bernoulli polynomials have in relation with the distribution $p_{n}$ or the error $\delta_n$. It is not also clear from this recursion, even with \eqref{eq:SymmOfQ}, that $\delta_n$ has the symmetry described in \eqref{eq:deltaSymm}.
			
			\item Considering $\delta_n$ as a discrete analogue of $b_n$, one can think of the equation $\delta_{n+1} = \delta_{n} \ P_{n+1}$ as a discrete analogue of the integral formula \eqref{eq:intchar}.
			
			\item In addition to the dynamics of the Markov chain $I = (I_1, I_2, \dots)$, we can get obtain the joint distribution of $(I_n, D_n)$ recursively in the same way. The key observation is that at step $n$, having obtained $I_n$ from the Bernoulli clock scheme and inserting the two new values $n+1$ in the clock, we either increment $D_n$ by $1$ to get $D_{n+1}$ if both values are inserted prior to $I_n$ or the number of laps is not incremented i.e. $D_{n+1} = D_n$ if one of the two values is inserted after $I_n$. We then obtain the following recursion for $\#(n; i, d)$: 
			
			\begin{enumerate}
				\item[1)] $\#(1  ; 1,   0) = 1$
				\item[2)] $\#(n+1; i, d) = \sum\limits_{1 \leq  x < h}  \#(n; i, x) \ \#_{n+1}(x,d) + \sum\limits_{h \leq x \leq 2n}   \#(n; i-1, x) \ \#_{n+1}(x,d)$.
			\end{enumerate}
			So one can get the joint distribution of $(I_n,D_n)$ recursively with
			\[
			\P(I_{n} = i, \ D_n = d) = \frac{\#(n; i, d )}{2^{-n}(2n)!}.
			\] 
		\end{enumerate}
		
	\end{remark}

	\section{Generalized Bernoulli clock} \label{Sec:Generalized_Bernoulli_clock}
	
	Let $n \geq 1$, $m_1, \dots, m_n \geq 1$ be positive integers and $M = m_1 + \dots + m_n$. 
	Let $\tau_n = \tau(m_1, \dots, m_n)$ be a random permutation uniformly distributed among the $M!/(m_1 ! \dots m_n !)$ permutations of the multiset $1^{m_1}2^{m_2} \dots n^{m_n}$.  
	Let us denote by $1 \leq I_1 \leq M$ the index of the first $1$ in the sequence $\tau_n$. 
	Continuing from this index $I_1$, let $I_2$ be the index of the first $2$ we encounter (circling back if necessary) and continuing in this manner we get random indices $(I_1, I_2, \dots, I_n)$. 
	Let us denote by $D_n = D(m_1, \dots, m_n)$ the number of times we circled around the sequence $\tau_n$ in this process,
	that is the number of descents in the random sequence $(I_1, I_2, \dots, I_n)$, as in \eqref{eq:dndef}.
	
	For the continuous model, mark the circle $\T = \R/\Z \cong [0,1)$  with $M$ i.i.d uniform on $[0,1]$ random variables $U_{1}^{(1)}, \dots, U_{1}^{(m_1)}$ , \dots , $U_{n}^{(1)}, \dots, U_{n}^{(m_n)}$ and let $U_{1:M} < \dots < U_{M:M}$ be their order statistics. Starting from $0$ we walk around the clock until we encounter the first of the variables $U_{1}^{(i)}$ at some random index $I_1$. We continue from the random index $I_1$ until we encounter the first of the variables $U_{2}^{(i)}$ (circling back if necessary) and continue like this until we encounter the first of the variables $U_{n}^{(i)}$. We then obtain the random sequence $(I_1, I_2, \dots, I_n)$ and $D_{n}$ is the number of times we circled around the clock. Finally, let us denote by $(X_1, \dots, X_n)$ the lengths (clock-wise) of the segments $[U_{I_1:M}, U_{I_2:M}]$, $\cdots$, $[U_{I_{n-1}:M}, U_{I_{n}:M}]$, $  [U_{I_{n}:M}, U_{I_{1}:M}]$ on the clock.  The model described in \cref{Sec:BernoulliClock} is the particular instance of this model where $m_1 = \dots = m_n = 2$.

	\begin{remark}
		When there is no risk of confusion, we shall suppress the parameters $m_1, \dots, m_n$ to simplify the notation.
	\end{remark}

	\begin{proposition} \label{prop:BernoulliClockGenCase}
		
		The following hold
		\begin{enumerate}
			\item The random lengths $X_1, X_2, \dots , X_n$ are independent random variables and $X_i$ has distribution beta($1,m_i$) for each $ 1 \leq i \leq n$.
			\item The random sequence of indices $(I_1, I_2, \dots, I_n)$ is independent of the order statistics $(U_{1:M} < \dots < U_{M:M})$.
		\end{enumerate}
		
	\end{proposition}
	\begin{proof}
		Notice that $X_1 = \min(U_1^{(1)}, \dots, U_1^{(m_1)})$ is a beta($1,m_1$) random variable. Also, since $U_2^{(1)}, \dots, U_2^{(m_2)}$ are i.i.d uniform and are independent of the position of $X_1$ on the circle, the variables
		$U_2^{(i)} - X_1 \mod \Z \in [0,1)$ are still i.i.d uniform so $X_2$ is also beta($1,m_2$) and independent of $X_1$. Running the same argument repeatedly we deduce that the variables $X_1, X_2 \dots, X_n$ are independent with $X_i \sim \mathrm{beta}(1,m_i)$. Also, the random index $I_n$ at which the process stops depends only on the relative positions of the variables $U_{1}^{(1)}, \dots, U_{1}^{(m_1)}$, \dots , $U_{n}^{(1)}, \dots, U_{n}^{(m_n)}$ i.e. $I_n$ is fully determined by the random permutation of $\{1,\dots, M\}$ induced by the $M$ i.i.d uniforms. We then deduce that $I_n$ is independent of the order statistics $(U_{1:M} < U_{2:M} < \dots < U_{M:M})$.
	\end{proof}
	
	The number $D_{n}$ of turns around the clock  can also be expressed as follows		
	\begin{equation}\label{eq:NumberOfLapSumoFBetas}
		D_{n} = \lfloor S_{n} \rfloor, \quad \text{ where } S_{n} \coloneqq X_1 + \dots + X_n.
	\end{equation}
	
	Let us denote by $L_n = L(m_1, m_2, \dots, m_n)$ the length of the longest continuous increasing subsequence of $\tau_n$ starting with $1$; that is the largest integer $1\leq \ell \leq n$ such that $$  1,2,3,\dots, \ell \quad \text{ is a subsequence of }\tau_n.  $$
	\begin{example}
		Suppose $n = 4$ and $(m_1,m_2,m_3,m_4) = (2,3,2,4)$ and consider the permutation $\tau_n = ( {1},4,4,1,4, {2 },4,  {3},3,2,2)$. The the longest increasing continuous subsequence of $\tau_n$ starting from $1$ (the boldfaced  subsequence) has length $L_4 = 3$ in this case.
	\end{example}

	For an infinite sequence $m = (m_1,m_2, \cdots )$ of positive integers, notice that we can construct the sequences of variables $L_n = L(m_1,\dots, m_n)$ , $D_n = D(m_1,\dots, m_n)$ and $I_n = I(m_1, \dots, m_n)$ on a common probability space. This is done by marking an additional $m_n$ i.i.d uniform positions on the circle $\T$ at each step $n$. Notice then that $(L_n = L(m_1, \dots, m_n))_{n \geq 1}$ is an increasing sequence of random variables so we define 
	\[
	L_{\infty} \coloneqq \lim\limits_{n \to \infty} L_n \quad \text{and} \quad  \Lcal_m  \coloneqq \E[L_\infty].
	\]
	\begin{proposition} \label{prop:LengthsAndLapNumbers}
		We have the following
		\[
		L_{n}  = \sum_{ k = 0 }^{n} 1[S_k \leq 1] \quad \text{and} \quad L_{\infty} = \sum_{k = 0}^{\infty} 1[S_k \leq 1].
		\]
		In particular, we have $(L_n = n) = (D_n = 0)$ and for $n \geq k$ we have
		\[
		( L(m_1, \dots , m_n) \geq k )	 = (L(m_1, \dots, m_k) = k).
		\]
	\end{proposition}   
	\begin{proof}
		The length $L_n$ of the longest sequence of the form $1 \dots \ell$ is the maximal integer $\ell$ such that $S_\ell \leq 1$, i.e. the maximal $l$ such that the random walk $(S_k)_{k\geq 0}$ does not shoot over $1$. Then we deduce that indeed
		\[
		L_{n}  = \sum_{ k = 0 }^{n} 1[S_k \leq 1].
		\]
		The rest of the statements follow immediately from this equation.
	\end{proof}
	
	\begin{corollary}
		For $k \leq n$ we have
		\[
		\P(L_n \geq k) = \P(S_{k} \leq 1).
		\]
	\end{corollary}
	\begin{proof}
		Follows immediately from \cref{prop:LengthsAndLapNumbers}.
	\end{proof}
	
	\begin{remark}
		When $m_1 = m_2 = \dots m_n = 1$, the random variable $S_n$ is the sum of $n$ i.i.d uniform random variables on $[0,1]$ and the fractional part $S_n^\circ$ has uniform distribution on $\T$. The index $I_n$ has uniform distribution in $\{1,\dots, n\}$ and the distribution of the number of descents
		$$P(D_n = k) = \frac{A_{n,k}}{n!}, \quad (0 \leq k \leq n-1)$$ 
		is given by the Eulerian numbers $A_{n,k}$, see \cite[Section 1.4]{StanleyVol1}.
	\end{remark}
	
	Horton and Kurn \cite[Theorem and Corollary (c)]{HortonKurn} gives a formula for the number of permutations $\tau$ of the multiset $1^{m_1}2^{m_2}\dots n^{m_n}$ for which $L_n = n$; that is a formula for
	\[
	\frac{M!}{m_1! \cdots m_n!} \ \P(L_n = n).
	\]
	We shall interpret this formula in our context and rederive  it from a probabilistic perspective.

	\begin{theorem}
		The number of permutations $\tau_n$ of the multiset $1^{m_1} 2 ^{m_2} \dots n^{m_n}$ that contain the sequence $(1,2,\cdots,n)$ is given by
		
		\begin{equation}\label{eq:numberOfPermutations}
			\frac{M!}{m_1! \cdots m_n!} \P(L_n = n) =  (-1)^M \sum_{j = 0}^{M}  \binom{M}{j} \frac{c_{j}}{j!},
		\end{equation}
		where
		\[
		c_j =  (-1)^n [\theta^j] \prod_{i = 1}^{n} E_{m_i - 1}(- \theta),
		\]
		with $[x^n] f(x)$ denoting the coefficient of $x^n$ in the power series expansion of $f$.
	\end{theorem}
	\begin{proof}
		Similarly to our discussion in \cref{Sec:Prob_Interpretations}, we can obtain an expression for $\P(S_{n} \leq x)$ by inverting the Laplace transform of $S_{n}$. First recall that the Laplace transform of $X_i \sim \mathrm{beta}(1,m_i)$ is
		\[
		\varphi_{X_i}(\theta) = \E[ e^{- \theta X_i }] =   (-1)^{m_i} \frac{m_i!}{\theta^{m_i}} \left( e^{-\theta} - E_{m_i - 1}(- \theta) \right),
		\]
		where $E_k(x)$ denotes the exponential polynomial $E_k(x) = \sum\limits_{i = 0}^{k} x^i/i!$. So the Laplace transform of $S_{n}$ is then given by
		\[
		\varphi_{S_{n}}(\theta) = \frac{(-1)^M \prod\limits_{i=1}^{n} m_i !} {\theta^M} \prod_{i = 1}^{n} \left( e^{-\theta} - E_{m_i - 1}(-\theta)\right).
		\]
		Using \eqref{eq:masterLaplaceTransform} to invert this Laplace transform, we get
		\[
		\P(S_{n} \leq x) =   (-1)^M \left(\prod\limits_{i=1}^{n} m_i ! \right) \sum_{k, j \geq 0} \alpha_{k,j} \frac{(x-k)^{M-j}_+ }{(M-j)!},
		\]
		where $\alpha_{k,j}$ is the coefficient of $\theta^j X^k$ in the polynomial $\prod_{i=1}^n (X - E_{m_i - 1}(-\theta))$. So we deduce that
		\[
		\P(L_n = n) = \P(S_n \leq 1) =  (-1)^M \left(\prod\limits_{i=1}^{n} m_i ! \right) \sum_{j = 0}^{M}  \frac{c_{j}}{(M-j)!},
		\]
		with
		\[
		c_j = \alpha_{0,j} =  (-1)^n [\theta^j] \left( \prod_{i = 1}^{n} E_{m_i - 1}(- \theta) \right).
		\]
		Multiplying by  $ M! / (m_1! \cdots m_n!)$ we get the formula \eqref{eq:numberOfPermutations}.
	\end{proof}

	We suppose from now on that $m \coloneqq m_1 = m_2 = \dots \geq 1$. Let $\Lcal_{n,m}$ and $\Lcal_{m}$ denote the expectation of $L_n$ and $L_\infty$; that is
	\[
	\Lcal_{n,m} \coloneqq \E[L_n] \quad \text{and} \quad \Lcal_{m} \coloneqq \lim_{n \to \infty} \Lcal_{n,m} = \E[L_{\infty}].
	\]
	In \cite{CliftonEtAl}, the authors present a fine asymptotic study of $\Lcal_{m}$ as $m \to \infty$. In this paper, we provide a pleasant probabilistic framework in which the discussion \cite{CliftonEtAl} fits rather naturally. 
	
	\medskip    
	
	Let $(N(t), t \geq 0)$ be the renewal process with beta($1,m$)-distributed i.i.d jumps $X_i$ i.e.
	\[
	N(t)  = \sum_{n \geq 1}^{\infty} 1[S_n \leq t].
	\]
	Notice that, by virtue of \cref{prop:LengthsAndLapNumbers}, the  variable $N(1) = L_\infty - 1 $ is the number of renewals of $N$ in $[0,1]$. Let $M(t) \coloneqq \E[N(t)]$ denote the mean of $N(t)$. By first step analysis, $M(t)$ satisfies the following equation for $t \in [0,1]$:
	
	\begin{align}\label{eq:funcEquationForM}
		M(t) &=  \P(X_1 \leq t) + m \int_{0}^{t} M(t-x) (1-x)^{m-1} dx, \\ 
		&= \P(X_1 \leq t) + m \int_{0}^{t} M(x) (1 - t + x)^{m-1} dx. \nonumber
	\end{align}
	
	From \eqref{eq:funcEquationForM} we can deduce that $M$ satisfies the following differential equation
	\begin{equation}\label{eq:DiffEquationForM}
		1 + \sum_{k = 0}^{m} \frac{(-1)^{k}}{k!} M^{(k)}(t)  = 0.
	\end{equation}
	
	\begin{theorem}\label{thm:ExpressionForMean}
		Let $\alpha_1, \dots, \alpha_m$ be the $m$ distinct complex roots of the exponential polynomial $E_m(x) = \sum_{k = 0}^{m} x^k/k!$. Then the mean function $M(t)$ is given by
		\begin{equation}\label{eq:FormulaForM}
			M(t) = - 1 - \sum_{k= 1}^{m} \alpha_k^{-1} e^{-\alpha_k t}. 
		\end{equation}					
	\end{theorem}
	
	Before we prove \cref{thm:ExpressionForMean}, we first recall a couple of intermediate results.
	
	\begin{lemma}\label{lem:Integ}
		Let $z$ be a non-zero complex number. Then, for any positive integer $n$ and $t \in [0,1]$, we have the following:
		\[
		\int_{0}^{t} e^{z x} (1-x)^{n} dx = n! \sum_{j = 0}^{n}  \frac{ e^{zt}(1-t)^j - 1}{j!} z^{j-n-1}.
		\]
	\end{lemma}
	\begin{proof}
		Follows immediately by induction on $n$ and integration by parts.
	\end{proof}
	The following lemma is an adaptation of \cite[Theorem 7]{Zemyan05}.
	\begin{lemma}\label{lem:Zemyan05}
		Let $\alpha_1, \dots, \alpha_m$ be the $m$ distinct complex zeros of $E_m(x)$. Then we have the following
		\[
		\sum_{k = 1}^{m} \alpha_k^{-j} = \begin{cases} -1, \quad &\text{if } j = 1,\\
			0, \quad &\text{if } 2 \leq j \leq m,\\
			1/m!, \quad &\text{if } j = m+1.
		\end{cases}
		\]
	\end{lemma}

	\begin{proof}[Proof of \cref{thm:ExpressionForMean}]
		The mean function $M(t)$ satisfies \eqref{eq:DiffEquationForM}. The latter is an order $m$ ODE with constant coefficients and its characteristic polynomial is $E_m(-x)$ whose roots are $-\alpha_1, \dots , - \alpha_m$. So the solution is of the form
		\[
		M(t) = -1 + \sum_{k = 1}^{m} \beta_k e^{-\alpha_k t}.
		\]
		Setting $\beta_k = - \alpha_k^{-1}$ for $1 \leq k \leq m$, it suffices to show that $M(t)$ satisfies \eqref{eq:funcEquationForM}. To that end notice that, thanks to \cref{lem:Integ}, we have
		\begin{align*}
			&\P(X_1 \leq t) +  m \int_{0}^{t} M(t-x) (1-x)^{m-1} dx\\
			&= m \int_{0}^{t} (1+M(t-x)) (1-x)^{m-1} dx\\
			&= - \sum_{k = 1}^{m}  m \alpha_k^{-1} \int_{0}^{t} e^{-\alpha_k (t-x)} (1-x)^{m-1} dx \\
			&= - \sum_{k = 1}^{m}  m \alpha_k^{-1} e^{-\alpha_k t} \int_{0}^{t} e^{\alpha_k x} (1-x)^{m-1} dx \\
			&=  \sum_{k = 1}^{m}  m \alpha_k^{-1} e^{-\alpha_k t}   (m-1)! \sum_{j = 0}^{m-1}  \frac{1 - e^{\alpha_k t}(1-t)^j}{j!} \alpha_k^{j - m}\\ 
			&=  m! \sum_{k = 1}^{m} \sum_{j = 0}^{m-1} \frac{e^{-\alpha_k t} - (1-t)^j}{j!} \alpha_k^{j-m-1}.
		\end{align*}
		Now notice that, thanks to \cref{lem:Zemyan05}, we have
		\[
		\sum_{j=0}^{m-1} \sum_{k=1}^{m} \frac{(1-t)^j}{j!} \alpha_k^{j-m-1} =  \sum_{k=1}^{m} \alpha_k^{-m-1} = \frac{1}{m!}.
		\]
		We also have
		\[
		\sum_{k=1}^{m} \sum_{j=0}^{m-1} \frac{e^{-\alpha_k t}}{j!} \alpha_k^{j-m-1} =  \sum_{k=1}^{m} \alpha_k^{-m-1} e^{-\alpha_k t} \sum_{j=0}^{m-1}  \frac{\alpha^j}{j!}  = - \frac{1}{m!} \sum_{k=1}^{m} \alpha_k^{-1} e^{-\alpha_k t}.
		\]
		The last equation follows from the fact that $\alpha_k$ is a zero of $E_m(x) = \sum_{j = 0}^{m} x^j / j!$.
		So combining the last two equations with the previous one, we get
		\[
		\P(X_1 \leq t) +  m \int_{0}^{t} M(t-x) (1-x)^{m-1} dx = - 1 - \sum_{k=1}^{m} \alpha_k^{-1} e^{-\alpha_k t} = M(t). \qedhere
		\]
	\end{proof}
	
	\begin{corollary}[Theorem 1.1-(a) in \cite{CliftonEtAl}]
		The expectation $\Lcal_{m}$ is given by 
		\[
		\Lcal_{m} = \sum_{k = 1}^{m} - \alpha_k^{-1} e^{ - \alpha_k}.
		\]
		In particular we have
		\[
		\Lcal_2 = e(\cos(1) + \sin(1)).
		\]
		
	\end{corollary}
	\begin{proof}
		Since $L_\infty = 1 + N(1)$, we deduce that $\Lcal_{m} = 1 + M(1)$ and the result follows immediately from \cref{thm:ExpressionForMean}.
	\end{proof}
	
	\begin{remark}
		Note that derivatives of $M$ at $0$ are the moments of the roots $\alpha_1, \dots, \alpha_m$ i.e.
		\[ \mu(j,m) \coloneqq  \sum_{k=1}^{m} \alpha_k^{j} = (-1)^j M^{(j+1)}(0), \quad \text{for }  j \geq 0. \]
		The functional equation \eqref{eq:funcEquationForM} then gives a recursion that these moments satisfy:
		\[ \mu(-1,m) = 0 \quad \text{and} \quad   \mu(j,m) =   (m)_{j+1} - \sum_{i = 0}^{j-1} \ (m)_{i+1}  \mu(j-i-1,m), \quad \text{for } j \geq 0. \]
		where $(X)_k = X(X-1)\cdots(X-k+1)$ is the $k$-th falling factorial polynomial. These moments are polynomials $\mu(j,\cdot)$ in $m$ and it would be interesting to give an expression for $\mu(j,X)$ and study its properties as suggested in \cite{Zemyan05}.
	\end{remark}

	To conclude this section, we give a positive answer to Conjectures 4.1 and 4.2 of \cite{CliftonEtAl}. For any integer $m \geq 1$, let $X_1^{(m)}, X_2^{(m)}, \dots $ be a sequence of i.i.d random variables with beta($1,m$) distribution and denote by $L_{n,m}$ and $L_{\infty, m}$ the following random variables
	\[
	L_{n,m} = \sum_{k = 1}^{n} 1\left[S^{(m)}_k \leq 1\right]  \quad \text{and} \quad L_{\infty, m} = \sum_{k = 1}^{\infty} 1\left[S^{(m)}_k \leq 1\right],
	\]
	with 
	\[
	S^{(m)}_n  = X_1^{(m)} + \cdots + X_n^{(m)}, \quad \text{for } n \geq 1.
	\]
	
	\begin{proposition}\label{prop:CLT}
		The random variable $(L_{\infty, m} - m) / \sqrt{m}$ converges in distribution to a Gaussian measure with mean $0$ and variance $1$.
	\end{proposition}   
	\begin{proof}
		For $m \geq 1$ and $x \in \R$ let $u(x,m) \coloneqq \lfloor m + x \sqrt{m} \rfloor $. We then have
		\begin{align*}
			\P\left( \frac{L_{\infty, m} - m}{\sqrt{m}}  \leq x \right) 
			&= \P\left(L_{\infty, m} \leq m + x\sqrt{m}\right)\\ 
			&= \P\left(L_{\infty, m} \leq u(x,m) \right)\\
			&= \P\left(S_{u(x,m) + 1} >  1 \right)\\
			&=\P\left( \frac{ m S_{u(x,m) + 1} - u(x,m) }{\sqrt{u(x,m)}}  > \frac{ m - u(x,m)}{\sqrt{u(x,m)}}  \right).
		\end{align*}
		Denote by $(Y_{k,m})$ the array defined as follows:
		\[
		Y_{k,m} = \frac{1}{\sqrt{m}} (m X_k^{(m)} - 1), \quad \text{for } k, m \geq 1.
		\]
		We then have $\E[Y_{k,m}] = 0$ and this array satisfies the conditions for the Lindeberg-Feller theorem \cite[Theorem 3.4.10]{Durett}, see \cref{appendix:CLT}. Applying this theorem yields        
		\[
		Y_{m,1} + \dots + Y_{m, m}  \xrightarrow[m \uparrow \infty]{} \mathcal{N}(0,1),
		\]
		but since $m - u(x,m) \sim x \sqrt{m}$ as $m \uparrow \infty$ we also deduce that 
		\[
		Y_{m,1} + \dots + Y_{m, u(x,m)}  \xrightarrow[m \uparrow \infty]{} \mathcal{N}(0,1).
		\]
		To conclude, notice that:
		\[
		\frac{ m S_{u(x,m) + 1} - u(x,m) }{\sqrt{u(x,m)}} = \sqrt{ \frac{m}{u(x,m)} } \left( Y_{m,1} + \dots + Y_{m, u(x,m)} \right)  \to \mathcal{N}(0,1) \quad \text{as } m \uparrow \infty.
		\]
		and  
		\[
		\frac{m - u(x,m)}{\sqrt{u(x,m)}}  \to -x \quad \text{as } m \uparrow \infty.
		\]
		So we deduce:
		\[
		\P\left( \frac{L_{\infty, m} - m}{\sqrt{m}}  \leq x \right)  \to  \int_{-x}^{\infty} \frac{1}{\sqrt{2\pi}} e^{-u^2/2} du =   \int_{- \infty}^{x} \frac{1}{\sqrt{2\pi}} e^{-u^2/2} du, \quad \text{as } m \uparrow \infty \qedhere
		\]
	\end{proof}

	\section{Wrapping probability distributions on the circle}\label{Sec:convOfSemigroups}
	In the decomposition \eqref{Xdec} for an exponentially distributed
	$X = \gamma_1/\lambda$ with parameter $\lambda > 0$; that is
	\[
	\P(X > t) = e^{- \lambda t}, \quad \text{for } t \geq 0,
	\]
	the Eulerian generating function \eqref{eq:eulersym} is the probability density of the fractional part $(\gamma_1/\lambda)^\circ$ at $u \in [0,1)$.
	In this probabilistic representation of Euler's exponential generating function \eqref{eq:Euler_Gen_Func},
	the factorially normalized Bernoulli polynomials $b_n(u)$ for $n>0$ are the densities at $u \in [0,1)$ of  a sequence of 
	signed measures on $[0,1)$, each with total mass $0$, which when weighted by $(-\lambda)^n$ and summed over $n >0$ give the difference between the 
	probability density of $(\gamma_1/\lambda)^\circ$ and the uniform probability density $b_0(u) \equiv 1$ for $u \in [0,1)$.

	For a positive integer $r$ and a positive real number $\lambda$, let $f_{r,\lambda}$ denote the probability density of the gamma($r,\lambda$) distribution:
	\[
	f_{\gamma_{r,\lambda}}(x) = \frac{\lambda^r}{\Gamma(r)} e^{-\lambda x} x^{r - 1} 1_{x > 0}, \quad x \in \R.
	\]
	It is well known that $f_{r,\lambda}$ is the $r$-fold convolution of $f_{1,\lambda}$ on the real line i.e. $f_{\gamma_{r,\lambda}} = (f_{\gamma_{1,\lambda}})^{ \ast r}$.

	Let $\gamma_{r,\lambda}$ be a random variable with distribution gamma($r,\lambda$) and let us denote by $\gamma_{r,\lambda}^\circ$ the random variable $\gamma_{r,\lambda} \mod \Z $ on the circle $\T$. The probability density of $\gamma_{r,\lambda}^\circ$ on $\T = [0,1)$ is given for $0 \le u < 1$ by
	\begin{align}\label{wrapped_gamma}
		f_{\gamma_{r,\lambda}}^\circ(u) = \sum_{m \in \Z}^{} f_{\gamma_{r,\lambda}}(u + m) &= \frac{\lambda^r}{\Gamma(r)} e^{-\lambda u}  \sum_{m = 0 }^{\infty}  (u+m)^{r-1} e^{-\lambda m}  \\
		&= \frac{\lambda^r}{\Gamma(r)} e^{-\lambda u } \Phi(e^{-\lambda}, 1 -r, u), \label{eq:wrappedGammaLerchHurwitz}
	\end{align}
	where $\Phi$ is the Hurwitz-Lerch zeta function $\Phi(z,s,u) = \sum_{m \geq 0} \frac{z^m}{(u + m)^s}$. In particular, for $r = 1$ the probability density of 
	$\gamma_{1,\lambda}^\circ$, the fractional part of an exponential variable with mean $1/\lambda$, at $u \in [0,1)$, is 
	\[
	f_{\gamma_{1,\lambda}}^\circ(u) = \frac{\lambda e^{\lambda(1-u)}}{e^{\lambda } - 1} = B(1-u,\lambda) = 1 + \sum_{n = 1}^{\infty} 
	b_n(1-u) \lambda^n
	\]
	where $B(x,\lambda)$, evaluated here for $x = 1-u$, is the generating function in \eqref{eq:Euler_Gen_Func}. Combined with the reflection symmetry \eqref{eq:refl}, this shows that the probability density of  $\gamma_{1,\lambda}^\circ$ can be expanded in Bernoulli polynomials as:
	\begin{equation} \label{eq:wrapped_exp}
		f_{\gamma_{1,\lambda}}^\circ(u) = 1 + \sum_{n = 1}^{\infty} (-1)^n b_n(u) \lambda^n \qquad (0 \leq u < 1 ).
	\end{equation}
	The following proposition generalizes this result to all integers $r \geq 1$.

	The expansion \eqref{eq:expansion_gamma} can be read from \eqref{eq:wrappedGammaLerchHurwitz} and formula (11) on page 30 of \cite{MR0058756}. The consequent interpretation \eqref{eq:norlund} of $b_r(u)$ for $r >0$, as the density of a signed measure describing how the probability density $f_{\gamma_{r,\lambda}}^\circ(u)$ approaches the uniform density $1$ as $\lambda \downarrow 0$, dates back to the work of N\"orlund \cite[p. 53]{norlund_vorlesungen_1924}, who gave an entirely analytical account of this result. 
	See also \cite{Coelho} for further study of the wrapped gamma and related probability distributions, and \cite{DSVIdentities} for various identities related to \eqref{eq:expansion_gamma}.

	\begin{proposition}[Wrapped gamma distribution] \label{prop:wrappedGamma}
		For each $r = 1,2,3, \ldots$ the wrapped gamma density admits the following expansion:
		\begin{equation} \label{eq:expansion_gamma}
			f_{\gamma_{r,\lambda}}^\circ(u) = 1 + \sum_{n=r}^{\infty} (-1)^{n-r+1} \binom{n-1}{r-1} b_n(u) \lambda^n  \quad \text{ for }  0 < \lambda < 2 \pi 
		\end{equation}
		where the convergence is uniform in $u \in [0,1)$. In particular, as $\lambda \downarrow 0$ 
		\begin{equation} \label{eq:norlund}
			f_{\gamma_{r,\lambda}}^\circ(u) = 1 - \lambda^r b_r(u) + O(\lambda^{r+1}),  \quad \textrm{ uniformly in } u \in [0,1).
		\end{equation}
	\end{proposition}
	
	\begin{proof}
		Since $f_{\gamma_{r,\lambda}} = (f_{\gamma_{1,\lambda}})^{ \circst r}$ we deduce that $f_{\gamma_{r,\lambda}}^\circ = (f_{\gamma_{1,\lambda}}^\circ)^{\circst r}$. Then, combining \eqref{eq:wrapped_exp} and \cref{cor:Bernoulliconv} we deduce that
		\begin{align*}
			f_{\gamma_{r,\lambda}}^\circ(u) &= (\underbrace{ f_{\gamma_{1,\lambda}}^\circ \circst \dots \circst f_{\gamma_{1,\lambda}}^\circ}_{r \textrm{ factors}})(u)\\
			&= 1 + \sum_{k_1,\dots, k_r \geq 1} (-1) ^{k_1 + \dots  + k_r}\lambda ^{k_1 + \dots  + k_r} (b_{k_1} \circst \dots \circst b_{k_r})(u)\\
			&= 1 + \sum_{n = r}^{\infty} \sum_{ \substack{  k_1,\dots,k_r \geq 1 \\  k_1+\dots+ k_r = n } }^{} (-1)^n \lambda ^n  (-1)^{-r+1 }b_{n}(u)\\
			&= 1+ \sum_{n = r}^{\infty}  (-1)^{n-r+1} A_{r,n} \lambda ^n b_{n}(u),
		\end{align*}
		where $A_{r,n} = \binom{n-1}{r-1}$ is the number of $r$-tuples of positive integers that sum to $n$. Notice that all the sums we considered are summable uniformly in $u \in [0,1]$ since $\norm{b_n}_\infty = O((2\pi)^n) $ as $n \to \infty$, see \eqref{eq:SineFunctionApprox}.
	\end{proof}
	
	\begin{remark}
		The general problem of expanding a function on $\T$ as a sum of Bernoulli polynomials was first treated Jordan \cite[Section 85]{jordanFiniteDiff} and Mordell \cite{mordell_expansion_1966}. In our context, we think of the expansion of a function in Bernoulli polynomials as an analog of the Taylor expansion where we work with the convolutions $\circst$ instead of the usual multiplication of functions; i.e. we view expansions of the form
		\[
		f(x) =  a_0(f) + \sum_{n = 1}^{\infty} (-1)^{n-1} a_{n}(f) b_1^{\circst n} (x) = a_0(f) + \sum_{n = 1}^{\infty} a_{n}(f)  b_n(x),
		\]
		as an analogue of Taylor expansions
		\[
		f(x) = f(0) + \sum_{n=1}^{\infty} \frac{f^{(n)}(0)}{n!} x^n.
		\]
		As we have seen in this section, this point of view is especially fruitful when one wishes to convolve probability measures on $\T = [0,1)$. If $f$ is a $C^\infty$ function on $[0,1]$ satisfying some dominance condition (see \cite[Theorem 1]{mordell_expansion_1966}), the coefficient of $b_1^{\circst}(x)$ in the expansion of $f$ is given by
		\[
		(-1)^{n-1} a_n(f) =  (f^{(n-1)}(1) - f^{(n-1)}(0)), \quad \text{for } n \geq 0.
		\]
	\end{remark}    
	
	\appendix
	
	\section{An elementary combinatorial proof of \texorpdfstring{\cref{thm:Bernoulliconv}}{}}
	\label{appendix:combinatorial proof}
	
	As promised in \cref{rem:}, we give an elementary combinatorial proof of \cref{thm:Bernoulliconv} using generating functions. We first recall the following identity of the Bernoulli numbers $B_n$:
	
	\begin{equation}\label{eq:Bern_coeffs}
		B_{n} = \frac{-1}{n + 1} \sum_{k = 0}^{n-1} \binom{n+1}{k} B_{k}, \quad \text{for }n \geq 1.      
	\end{equation}

	\begin{proof}[Proof of \cref{thm:Bernoulliconv}]
		We proceed by induction. The first two polynomials $B_0(x)$ and ${B_1}(x)$ obviously satisfy $\cref{thm:Bernoulliconv}$. For $n \geq 1$, assume that $ {B_n}(x) = (-1)^{n-1} n ! \  \underbrace{ {B_1(x)}  \circst \dots  \circst  {B_1(x)}}_{n \ \mathrm{ factors}}$. We want to show that
		\[
		{B_{n+1}}(x) = - (n+1) {B_1(x)}  \circst  {B_n(x)}.
		\]
		For this, we use \cref{prop:convMonomRec} to compute ${B_1}  \circst  {B_n}$ as follows:
		\begin{align*}
			{x}  \circst  {B_n}(x)   &= x  \circst \sum_{k=0}^{n} \binom{n}{k} B_{n-k} x^{k} \\
			&= B_n \ x  \circst  {1} + \sum_{k = 1}^{n} \binom{n}{k} B_{n-k}  {x}  \circst  {x^k}\\
			&= \frac{B_n}{2} + \sum_{k = 1}^{n} \binom{n}{k} B_{n-k} \left( \frac{x - x^{k+1}}{k+1} + \frac{1}{(k+1)(k+2)}\right)\\
			&=\sum_{k = 1}^{n} \binom{n}{k} B_{n-k} \frac{x - x^{k+1}}{k+1} + \sum_{k = 0}^{n} \frac{\binom{n}{k} B_{n-k}}{(k+1)(k+2)},
		\end{align*} 
		and since $n \geq 1$ we have $1 \circst B_n(x) = 0$. Given that $B_1(x) = x - 1/2$, we have
		\begin{equation}\label{eq:xConvB_n}
			(n+1) {B_1(x)}  \circst  {B_n(x)}  = \sum_{k = 1}^{n} (n+1) \binom{n}{k} B_{n-k} \frac{x - x^{k+1}}{k+1} + \sum_{k = 0}^{n} \binom{n}{k} \frac{ (n+1) B_{n-k}}{(k+1)(k+2)}.
		\end{equation}
		We now expand the latter polynomial to match the expansion of $B_{n+1}(x) = \sum_{k=0}^{n+1} \binom{n+1}{k} B_{n + 1 - k} x^{k}$. From \eqref{eq:xConvB_n} we deduce that
		\[
		\Scale[0.97]{  (n+1) {B_1}(x)  \circst  {B_n}(x) = - \sum\limits_{k = 2}^{n+1} \binom{n+1}{k} B_{n+1-k} x^k  + \left( \sum\limits_{k = 2}^{n+1} \binom{n+1}{k} B_{n+1-k} \right) x    +  \sum_{k = 0}^{n} \binom{n}{k} \frac{ (n+1) B_{n-k}}{(k+1)(k+2)}.} 
		\]
		Notice that, thanks to the recursion \cref{eq:Bern_coeffs}, the coefficient of $x$ in the polynomial $(n+1)x \circst B_n(x)$ is
		\[
		\sum_{k = 2}^{n+1} \binom{n+1}{k} B_{n+1-k} = \sum_{k = 0}^{n-1} \binom{n+1}{k} B_{k} = - (n+1) B_n.
		\]
		So we deduce that
		\begin{align*}
			(n+1)( {B_1}  \circst  {B_n})(x)  &= - \sum_{k = 2}^{n+1} \binom{n+1}{k} B_{n+1-k} x^k  - (n+1) B_n x +  \sum_{k = 0}^{n} \binom{n}{k} \frac{ (n+1) B_{n-k}}{(k+1)(k+2)}\\
			& = - \sum_{k = 1}^{n+1} \binom{n+1}{k} B_{n+1-k} x^k  +  \sum_{k = 0}^{n} \binom{n}{k} \frac{ (n+1) B_{n-k}}{(k+1)(k+2)}.
		\end{align*}
		All that remains is to deal with the constant coefficient in \eqref{eq:xConvB_n}, and from \cref{lem:constantCoeff} we can see that the constant coefficient in the polynomial $(n+1)( {B_1}  \circst  {B_n})(x)$ is
		\[
		\sum_{k = 0}^{n} \binom{n}{k} \frac{ (n+1) B_{n-k}}{(k+1)(k+2)} = -(n+1) B_{n+1}.
		\]
		Hence, we obtain the desired equation
		\[	
		(n+1)( {B_1}  \circst  {B_n})(x) = - \sum_{k = 0}^{n+1} \binom{n+1}{k} B_{n+1-k} x^k = -  {B_{n+1}}(x),
		\]
		where the last equality is deduced from to \eqref{eq:Bern_coeffs}.
	\end{proof}
	
	\begin{lemma}\label{lem:constantCoeff}
		For any integer $n \geq 0$ the following equation holds:
		\[
		\sum_{k = 0}^{n} \frac{1}{(k+2)!} \frac{B_{n-k}}{(n-k)!}  =  - \frac{B_{n+1}}{(n+1)!}.
		\]
	\end{lemma}
	\begin{proof}
		The generating function of the sequence $\left(\frac{1}{(n+2)!}\right)_{n \geq 0}$ is the function
		\[
		g(z) \coloneqq \sum_{n = 0}^{\infty} \frac{z^n}{(n+2)!} = \frac{e^z - z - 1}{z^2},
		\]
		and the generating function of the sequence $\left(\frac{B_n}{n!}\right)_{n \geq 0}$ is $B(0,z) \coloneqq \sum_{n \geq 0} \frac{B_n}{n!} z^n = \frac{z}{e^z - 1}$. So the generating function of the convolution of the two sequences is 
		\[
		h(z) \coloneqq g(z)B(0, z) = \sum_{n=0}^{\infty} \left( \sum_{k = 0}^{n} \frac{1}{(k+2)!} \frac{B_{n-k}}{(n-k)!} \right) z^n =  \frac{e^z  - z - 1}{ z (e^z - 1)}.
		\]
		Now, the generating function of the sequence $\left ( \frac{B_{n+1}}{(n+1)!} \right)_{n\geq 0}$ is 
		\[
		f(z) \coloneqq \sum_{n=0}^{\infty}   \frac{B_{n+1}}{(n+1)!} z^{n} = \frac{B(z) - 1}{z} = \frac{z - e^z + 1}{z (e^z - 1)}.
		\]
		We deduce that $h(z) = - f(z)$ hence the desired result.
	\end{proof}

	\section{Complement to the proof of \texorpdfstring{\cref{prop:CLT}}{}}
	\label{appendix:CLT}
	Here we check that the array $Y_{k,m} = (m X_k^{(m)} - 1) / \sqrt{m} $ where $X_{1}^{(m)}, X_{2}^{(m)}, \dots $ is a sequence of i.i.d beta($1, m$) random variables, satisfies the conditions required in the Lindeberg-Feller theorem \cite[Theorem 3.4.10]{Durett}. For that we need to check the following:
	\begin{enumerate}
		\item $\sum_{k=1}^{m} \E[Y_{k,m}^2] \xrightarrow[]{m \to \infty} 1$.
		
		\item For any $\epsilon > 0$, we have $\sum_{k=1}^{m} \E[Y_{k,m}^2; |Y_{k,m}| > \epsilon] \xrightarrow[]{m \to \infty} 0$.
	\end{enumerate}
	For the first condition we have
	\[
	\sum_{k=1}^{m} \E[Y_{k,m}^2] = m^2 \mathrm{Var}(X_{k}^{(m)}) = \frac{m^3}{(m+1)^2 (m+2)}  \xrightarrow[m \to \infty]{} 1.
	\]
	For the second condition, fix $\epsilon > 0$ and note that the density of $Y_{k,m}$ is
	\[
	g_m(y) = \sqrt{m} \left( 1 -  \frac{\sqrt{m} \ y +1 }{m} \right)^{m-1},   \quad \text{for } -1/\sqrt{m} \leq y \leq (m-1)/\sqrt{m}.
	\]
	So for large enough $m$ we get
	\begin{align*}
		\E[Y_{k,m}^2; |Y_{k,m}| > \epsilon] 
		&= \int_{-1/\sqrt{m}}^{(m-1)/\sqrt{m}}  y^2 g_m(y) \ 1[|y|>\epsilon] dy \\
		&= \int_{\epsilon}^{(m-1)/\sqrt{m}}  y^2 g_m(y) dy \\
		&= \sqrt{m} \int_{\epsilon}^{(m-1)/\sqrt{m}}  y^2 \left( 1 -  \frac{\sqrt{m} \ y +1 }{m} \right)^{m-1}  dy. \\
	\end{align*}
	With the change of variable $z = (\sqrt{m}y + 1) / m$ we get
	\begin{align*}
		\E[Y_{k,m}^2; |Y_{k,m}| > \epsilon] 
		&= \int_{(\epsilon \sqrt{m} +1)/m}^{1} (mz - 1)^2 (1-z)^{m-1} dz\\
		&= \left(1 - \frac{\epsilon \sqrt{m} + 1}{m} \right)^m \frac{m( \epsilon^2 m (m+1) + 2m + 2 \epsilon \sqrt{m}(m-1) - 4 ) + 2 }{m(m+1)(m+2)}.
	\end{align*}
	So we deduce that 
	\begin{align*}
		\sum_{k=1}^{m} \E[Y_{k,m}^2; |Y_{k,m}| > \epsilon] &= \left(1 - \frac{\epsilon \sqrt{m} + 1}{m} \right)^m \frac{m( \epsilon^2 m (m+1) + 2m + 2 \epsilon \sqrt{m}(m-1) - 4 ) + 2 }{(m+1)(m+2)} \\
		& \substack{\simeq \\ m \to \infty} \ \epsilon^2 m e^{-\epsilon \sqrt{m}}.
	\end{align*}
	So we deduce that 
	\[
	\sum_{k=1}^{m} \E[Y_{k,m}^2; |Y_{k,m}| > \epsilon] \xrightarrow[m \to \infty]{} 0.
	\]
	
	\bibliographystyle{acm}
	\bibliography{refs}
	
\end{document}